\documentclass[12pt,twoside]{article}
\usepackage[english]{babel}
\usepackage[latin1]{inputenc}
\usepackage{amsmath}
\usepackage{amssymb,amsfonts}
\usepackage{graphicx}                   

\newcommand{\bR}{\mathbf{R}}

\newcommand{\bZ}{\mathbf{Z}}

\newcommand{\Bu}{\boldsymbol{u}}

\newcommand{\HYP}{\mathbb{H}^3}
\newcommand{\HYN}{\mathbb{H}^n}

\begin{document}
\pagestyle{myheadings}
\markboth{\centerline{Jen\H o Szirmai}}
{Packings with horo- and hyperballs \dots}
\title
{Packings with horo- and hyperballs generated by simple frustum orthoschemes}

\author{\normalsize{Jen\H o Szirmai} \\
\normalsize Budapest University of Technology and \\
\normalsize Economics Institute of Mathematics, \\
\normalsize Department of Geometry \\
\date{\normalsize{\today}}}

\maketitle


\begin{abstract}

In this paper we deal with the packings derived by horo- and hyperballs (briefly hyp-hor packings) in the $n$-dimensional hyperbolic spaces $\HYN$
($n=2,3$) which form a new class of the classical packing problems.

We construct in the $2-$ and $3-$dimensional hyperbolic spaces hyp-hor packings that
are generated by complete Coxeter tilings of degree $1$ i.e. the fundamental domains of these tilings are simple frustum orthoschemes
and we determine their densest packing configurations and their densities.

We prove that in the hyperbolic plane ($n=2$) the density of the above hyp-hor packings arbitrarily approximate
the universal upper bound of the hypercycle or horocycle packing density $\frac{3}{\pi}$ and
in $\HYP$ the optimal configuration belongs to the $[7,3,6]$ Coxeter tiling with density $\approx 0.83267$.

Moreover, we study the hyp-hor packings in
truncated orthosche\-mes $[p,3,6]$ $(6< p < 7, ~ p\in \bR)$ whose
density function is attained its maximum for a parameter which lies in the interval $[6.05,6.06]$
and the densities for parameters lying in this interval are larger that $\approx 0.85397$. That means that these
locally optimal hyp-hor configurations provide larger densities that the B\"or\"oczky-Florian density upper bound
$(\approx 0.85328)$ for ball and
horoball packings but these hyp-hor packing configurations can not be extended to the entirety of hyperbolic space $\mathbb{H}^3$.
\end{abstract}

\newtheorem{theorem}{Theorem}[section]
\newtheorem{corollary}[theorem]{Corollary}
\newtheorem{conjecture}{Conjecture}[section]
\newtheorem{lemma}[theorem]{Lemma}
\newtheorem{exmple}[theorem]{Example}
\newtheorem{defn}[theorem]{Definition}
\newtheorem{rmrk}[theorem]{Remark}
\newenvironment{definition}{\begin{defn}\normalfont}{\end{defn}}
\newenvironment{remark}{\begin{rmrk}\normalfont}{\end{rmrk}}
\newenvironment{example}{\begin{exmple}\normalfont}{\end{exmple}}
\newenvironment{acknowledgement}{Acknowledgement}


\section{Introduction}
The packing and covering problems with solely horo- or hyperballs (horo- or hypespheres) are intensively investigated
in earlier works in the $n$-dimensional $(n\ge2)$ hyperbolic space $\HYN$.
\begin{enumerate}
\item
{\bf On horoball packings}

In the $n$-dimensional hyperbolic space $\HYN$ there are $3$ kinds
of the "balls (spheres)" the balls (spheres), horoballs (horospheres) and hyperballs (hyperspheres).

The 2-dimensional case of circle and horocycle packings was settled by L.~Fejes T\'oth in \cite{FTL}.

The greatest possible density in hyperbolic space $\mathbb{H}^3$ is $\approx 0.85328$
which is not realized by packing regular balls. However, it is attained by a horoball packing of
$\overline{\mathbb{H}}^3$ where the ideal centers of horoballs lie on the
absolute figure of $\overline{\mathbb{H}}^3$. This ideal regular
simplex tiling is given with Coxeter-Schl\"afli symbol $[3,3,6]$ see e.g. \cite{Be}, \cite{B--F64}, \cite{B78} and \cite{G--K--K}.

In the previous paper \cite{KSz} we proved that the above known optimal ball packing arrangement in $\mathbb{H}^3$ is not unique.
We gave several new examples of horoball packing arrangements based on totally asymptotic Coxeter tilings that
yield the B\"or\"oczky--Florian packing density upper bound \cite{B--F64}.
Furthermore, by admitting horoballs of different types at each vertex of a totally asymptotic simplex and generalizing
the simplicial density function to $\mathbb{H}^n$ for $(n \ge 2)$ we find the B\"or\"oczky type density
upper bound is no longer valid for the fully asymptotic simplices in cases $n \ge 3$  \cite{Sz12}, \cite{Sz12-2}.
For example, the density of such optimal,
locally densest packing is $\approx 0.77038$ which is larger than the
analogous B\"or\"oczky type density upper bound of $\approx 0.73046$ for $\overline{\mathbb{H}}^4$.
However these ball packing configurations are only locally optimal and cannot be extended to the entirety of the hyperbolic spaces $\mathbb{H}^n$.

In the paper \cite{KSz14} we have continued our investigation of ball packings in hyperbolic 4-space
using horoball packings, allowing horoballs of different types.
We have shown seven counterexamples (which are realized by allowing one-, two-, or three horoball types)
to a conjecture of L. Fejes-T\'oth about the densest ball packings in hyperbolic $4$-space. The maximal density is $\approx 0.71645$

In \cite{Sz15} we proved that the optimal horoball density related to the hyperbolic 24 cell in $\mathbb{H}^4$ is $\approx 0.71645$ as well.
\item
{\bf{On hyperball packings}}

In \cite{Sz06-1} and \cite{Sz06-2} we have studied the regular prism tilings and the corresponding optimal hyperball packings in
$\mathbb{H}^n$ $(n=3,4)$ and in the paper  \cite{Sz13-3} we have extended the in former papers developed method
to 5-dimensional hyperbolic space and construct to each investigated Coxeter tiling a regular prism tiling,
have studied the corresponding optimal hyperball packings by congruent hyperballs,
moreover, we have determined their metric data and their densities. \newline
\indent In hyperbolic plane $\mathbb{H}^2$ the universal upper bound of the hypercycle packing density is $\frac{3}{\pi}$
proved by I.~Vermes in \cite{V79} and recently, (to the author's best knowledge) the candidates for the densest hyperball
(hypersphere) packings in the $3,4$ and $5$-dimensional hyperbolic space $\mathbb{H}^n$ are derived by the regular prism
tilings which are studied in papers \cite{Sz06-1}, \cite{Sz06-2} and \cite{Sz13-3}.\newline
\indent In $\mathbb{H}^2$ the universal lower bound of the hypercycle covering density is $\frac{\sqrt{12}}{\pi}$
determined by I.~Vermes in \cite{V81}. \newline
\indent In the paper \cite{Sz13-4} we have studied the $n$-dimensional $(n \ge 3)$ hyperbolic regular prism honeycombs
and the corresponding coverings by congruent hyperballs and we have determined their least dense covering densities.
Moreover, we have formulated a conjecture for the candidate of the least dense hyperball
covering by congruent hyperballs in the 3- and 5-dimensional hyperbolic space. \newline
In \cite{Sz14-2} we studied the problem of hyperball (hypersphere) packings in the
$3$-dimensional hyperbolic space.
We described to each saturated hyperball packing
a procedure to get a decomposition of the 3-dimensional hyperbolic space $\HYP$ into truncated
tetrahedra. Therefore, in order to get a density upper bound to hyperball packings it is sufficient to determine
the density upper bound of hyperball packings in truncated simplices.
We considered the hyperball packings in truncated simplices and proved that if the truncated tetrahedron is regular, then the
density of the densest packing is $\approx 0.86338$ which is larger than the B\"oröczky-Florian density upper bound,
however these hyperball packing configurations are only locally optimal and cannot be extended to the entirety of the hyperbolic spaces
$\mathbb{H}^3$.
\end{enumerate}
In this paper we deal with the packings with horo- and hyperballs (briefly hyp-hor packings) in the $n$-dimensional hyperbolic spaces $\HYN$
($n=2,3$) which form a new class of the classical packing problems.

We construct in the $2-$ and $3-$dimensional hyperbolic spaces hyp-hor packings that
are generated by complete Coxeter tilings of degree $1$ i.e. the fundamental domains of these tilings are
simple frustum orthoschemes with a principal vertex lying on the absolute quadric $Q$
and the other principal vertex is outer point.
We determine their densest packing configurations and their densities.
These Coxeter tilings exist in the $2-$, $3-$ and $5-$dimensional
hyperbolic spaces (see \cite{IH90}) and have given by their Coxeter-Schl\"afli graph in Fig.~1.
\begin{figure}[ht]
\centering
\includegraphics[width=12cm]{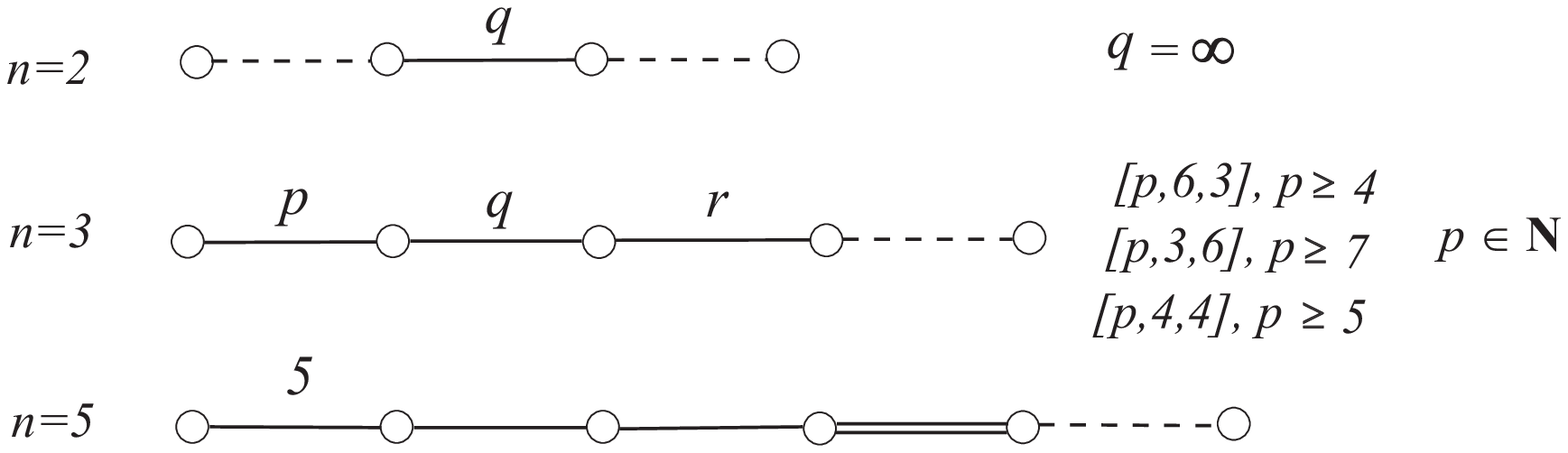}
\caption{Coxeter-Schl\"afli graph of Coxeter tilings of degree 1.}
\label{}
\end{figure}
We prove that in the hyperbolic plane $n=2$ the density of the above hyp-hor packings arbitrarily approximate
the universal upper bound of the hypercycle or horocycle packing density $\frac{3}{\pi}$ and
in $\HYP$ the optimal configuration belongs to the $[7,3,6]$ Coxeter tiling with density $\approx 0.83267$.

Moreover, we consider the hyp-hor packings in
truncated orthosche\-mes $[p,3,6]$ $(6< p < 7, ~ p\in \bR)$.
Its density function is attained its maximum for a parameter which lies in the interval $[6.05,6.06]$
and the densities for parameters lying in this interval are larger that $\approx 0.85397$. That means that these
locally optimal hyp-hor configurations provide larger densities that the B\"or\"oczky-Florian density upper bound
$(\approx 0.85328)$ for ball and
horoball packings but these hyp-hor packing configurations can not be extended to the entirety of hyperbolic space $\mathbb{H}^3$.
\section{The projective model and \\ the complete orthoschemes }
For $\mathbb{H}^n$ we use the projective model in the Lorentz space $\mathbb{E}^{1,n}$ of signature $(1,n)$,
i.e.~$\mathbb{E}^{1,n}$ denotes the real vector space $\mathbf{V}^{n+1}$ equipped with the bilinear
form of signature $(1,n)$
$
\langle ~ \mathbf{x},~\mathbf{y} \rangle = -x^0y^0+x^1y^1+ \dots + x^n y^n
$
where the non-zero vectors
$
\mathbf{x}=(x^0,x^1,\dots,x^n)\in\mathbf{V}^{n+1} \ \  \text{and} \ \ \mathbf{y}=(y^0,y^1,\dots,y^n)\in\mathbf{V}^{n+1},
$
are determined up to real factors, for representing points of $\mathcal{P}^n(\mathbb{R})$. Then $\mathbb{H}^n$ can be interpreted
as the interior of the quadric
$
Q=\{[\mathbf{x}]\in\mathcal{P}^n | \langle ~ \mathbf{x},~\mathbf{x} \rangle =0 \}=:\partial \mathbb{H}^n
$
in the real projective space $\mathcal{P}^n(\mathbf{V}^{n+1},
\mbox{\boldmath$V$}\!_{n+1})$.

The points of the boundary $\partial \mathbb{H}^n $ in $\mathcal{P}^n$
are called points at infinity of $\mathbb{H}^n $, the points lying outside $\partial \mathbb{H}^n $
are said to be outer points of $\mathbb{H}^n $ relative to $Q$. Let $P([\mathbf{x}]) \in \mathcal{P}^n$, a point
$[\mathbf{y}] \in \mathcal{P}^n$ is said to be conjugate to $[\mathbf{x}]$ relative to $Q$ if
$\langle ~ \mathbf{x},~\mathbf{y} \rangle =0$ holds. The set of all points which are conjugate to $P([\mathbf{x}])$
form a projective (polar) hyperplane
$
pol(P):=\{[\mathbf{y}]\in\mathcal{P}^n | \langle ~ \mathbf{x},~\mathbf{y} \rangle =0 \}.
$
Thus the quadric $Q$ induces a bijection
(linear polarity $\mathbf{V}^{n+1} \rightarrow
\mbox{\boldmath$V$}\!_{n+1})$)
from the points of $\mathcal{P}^n$
onto its hyperplanes.

The point $X [\bold{x}]$ and the hyperplane $\alpha [\mbox{\boldmath$a$}]$
are called incident if $\bold{x}\mbox{\boldmath$a$}=0$ ($\bold{x} \in \bold{V}^{n+1} \setminus \{\mathbf{0}\}, \ \mbox{\boldmath$a$} \in \mbox{\boldmath$V$}_{n+1}
\setminus \{\mbox{\boldmath$0$}\}$).
\begin{definition}
An orthoscheme $\mathcal{S}$ in $\mathbb{H}^n$ $(2\le n \in \mathbb{N})$ is a simplex bounded by $n+1$ hyperplanes $H^0,\dots,H^n$
such that
(see \cite{B--H, K89})
$
H^i \bot H^j, \  \text{for} \ j\ne i-1,i,i+1.
$
\end{definition}

{\it The orthoschemes of degree} $d$ in $\mathbb{H}^n$ are bounded by $n+d+1$ hyperplanes
$H^0,H^1,\dots,H^{n+d}$ such that $H^i \perp H^j$ for $j \ne i-1,~i,~i+1$, where, for $d=2$,
indices are taken modulo $n+3$. For a usual orthoscheme we denote the $(n+1)$-hyperface opposite to the vertex $A_i$
by $H^i$ $(0 \le i \le n)$. An orthoscheme $\mathcal{S}$ has $n$ dihedral angles which
are not right angles. Let $\alpha^{ij}$ denote the dihedral angle of $\mathcal{S}$
between the faces $H^i$ and $H^j$. Then we have
$
\alpha^{ij}=\frac{\pi}{2}, \ \ \text{if} \ \ 0 \le i < j -1 \le n.
$
The $n$ remaining dihedral angles $\alpha^{i,i+1}, \ (0 \le i \le n-1)$ are called the
essential angles of $\mathcal{S}$.
Geometrically, complete orthoschemes of degree $d$ can be described as follows:
\begin{enumerate}
\item
For $d=0$, they coincide with the class of classical orthoschemes introduced by
{{Schl\"afli}} (see Definitions 2.1).
The initial and final vertices, $A_0$ and $A_n$ of the orthogonal edge-path
$A_iA_{i+1},~ i=0,\dots,n-1$, are called principal vertices of the orthoscheme.
\item
A complete orthoscheme of degree $d=1$ can be interpreted as an
orthoscheme with one outer principal vertex, say $A_n$, which is truncated by
its polar plane $pol(A_n)$ (see Fig.~1 and 3). In this case the orthoscheme is called simply truncated with
outer vertex $A_n$.
\item
A complete orthoscheme of degree $d=2$ can be interpreted as an
orthoscheme with two outer principal vertices, $A_0,~A_n$, which is truncated by
its polar hyperplanes $pol(A_0)$ and $pol(A_n)$. In this case the orthoscheme is called doubly
truncated. We distinguish two different types of orthoschemes but I
will not enter into the details (see \cite{K91}).
\end{enumerate}
A $n$-dimensional tiling $\mathcal{P}$ (or solid tessellation, honeycomb) is an infinite set of
congruent polyhedra (polytopes) that fit together to fill all space $(\mathbb{H}^n~ (n \geqq 2))$ exactly once,
so that every face of each polyhedron (polytope) belongs to another polyhedron as well.
At present the cells are congruent orthoschemes.
A tiling with orthoschemes exists if and only if each dihedral angle of a tile is submultiple of $2\pi$
(in the hyperbolic plane the zero angle is also possible).

Another approach to describing tilings involves the analysis of their symmetry groups.
If $\mathcal{P}$ is such a simplex tiling, then any motion taking one cell into another maps the
entire tiling onto itself. The symmetry group of this tiling is denoted by
$Sym \mathcal{P}$.
Therefore the simplex is a fundamental domain of the group $Sym \mathcal{P}$ generated by reflections in its
$(n-1)$-dimensional hyperfaces.

The scheme of an orthoscheme $S$ is a weighted graph (characterizing $S \subset \mathbb{H}^n$
up to congruence) in which the nodes, numbered by $0,1,\dots,n$ correspond to the bounding hyperplanes of $\mathcal{S}$.
Two nodes are joined by an edge if the corresponding hyperplanes are not orthogonal.
\begin{figure}[ht]
\centering
\includegraphics[width=7cm]{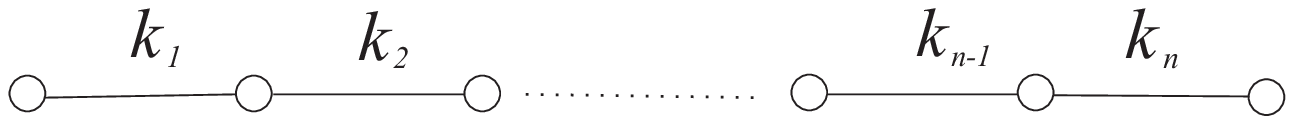}
\end{figure}
For the schemes of {\it complete Coxeter orthoschemes} $\mathcal{S} \subset \mathbb{H}^n$ we adopt the usual conventions and sometimes even use them in the Coxeter case: If two nodes are related by the weight $\cos{\frac{\pi}{p}}$
then they are joined by a ($p-2$)-fold line for $p=3,~4$ and by a single line marked $p$ for $p \geq 5$.
In the hyperbolic case if two bounding hyperplanes of $S$ are parallel, then the corresponding nodes
are joined by a line marked $\infty$. If they are divergent then their nodes are joined by a dotted line.

The ordered set $[k_1,\dots,k_{n-1},k_n] $ is said to be the
Coxeter-Schl$\ddot{a}$fli symbol of the simplex tiling $\mathcal{P}$ generated by $\mathcal{S}$.
To every scheme there is a corresponding
symmetric matrix $(c^{ij})$ of size $(n+1)\times(n+1)$ where $c^{ii}=1$ and, for $i \ne j\in \{0,1,2,\dots,n \}$,
$c^{ij}$ equals $-\cos{\frac{\pi}{k_{ij}}}$ with all angles between the facets $i$,$j$ of $\mathcal{S}$.

For example, $(c^{ij})$ below is the so called Coxeter-Schl\"afli matrix of the orthoscheme $S$ in
3-dimensional hyperbolic space $\mathbb{H}^3$ with
parameters (nodes) $k_1=p,k_2=q,k_3=r$ :
\[
(c^{ij}):=\begin{pmatrix}
1& -\cos{\frac{\pi}{p}}& 0 & 0 \\
-\cos{\frac{\pi}{p}} & 1 & -\cos{\frac{\pi}{q}}& 0 \\
0 & -\cos{\frac{\pi}{q}} & 1 & -\cos{\frac{\pi}{r}} \\
0 & 0 & -\cos{\frac{\pi}{r}} & 1 \\
\end{pmatrix}. \tag{2.4}
\]
\section{Basic notions and formulas}
\subsection{Coxeter tilings generated by simply frustum orthoschemes}
In general the complete Coxeter orthoschemes were classified by {{Im Hof}} in
\cite{IH85} by generalizing the method of {{Coxeter}} and {{B\"ohm}}, who
showed that they exist only for dimensions $\leq 9$. From this classification it follows, that the complete
orthoschemes of degree $d=1$ exist up to 5 dimensions.

In this paper we consider the orthoschemes of degree 1 where the initial vertex $A_0$ lies on the
absolute quadric $Q$. These orthoschemes and the corresponding Coxeter tilings exist in the $2$-, $3-$ and $5-$dimensional hyperbolic spaces and
are characterized by their Coxeter-Schl\"afli symbols and graphs (see Fig.~1).

In $n$-dimensional hyperbolic space $\mathbb{H}^n$ $(n \ge 2)$
it can be seen that if $\mathcal{S}=A_0A_1A_2 \dots A_n$ $P_0P_1P_2 \dots P_n$ is a complete
orthoscheme  with degree $d=1$ (a simply frustum orthoscheme) where $A_n$ is a outer vertex of
$\mathbb{H}^n$ then the points $P_0,P_1,P_2,\dots,P_{n-1}$ lie on the polar hyperplane $\pi$ of $A_n$ (see Fig.~2 in $\HYP$).
\begin{figure}[ht]
\centering
\includegraphics[width=6.5cm]{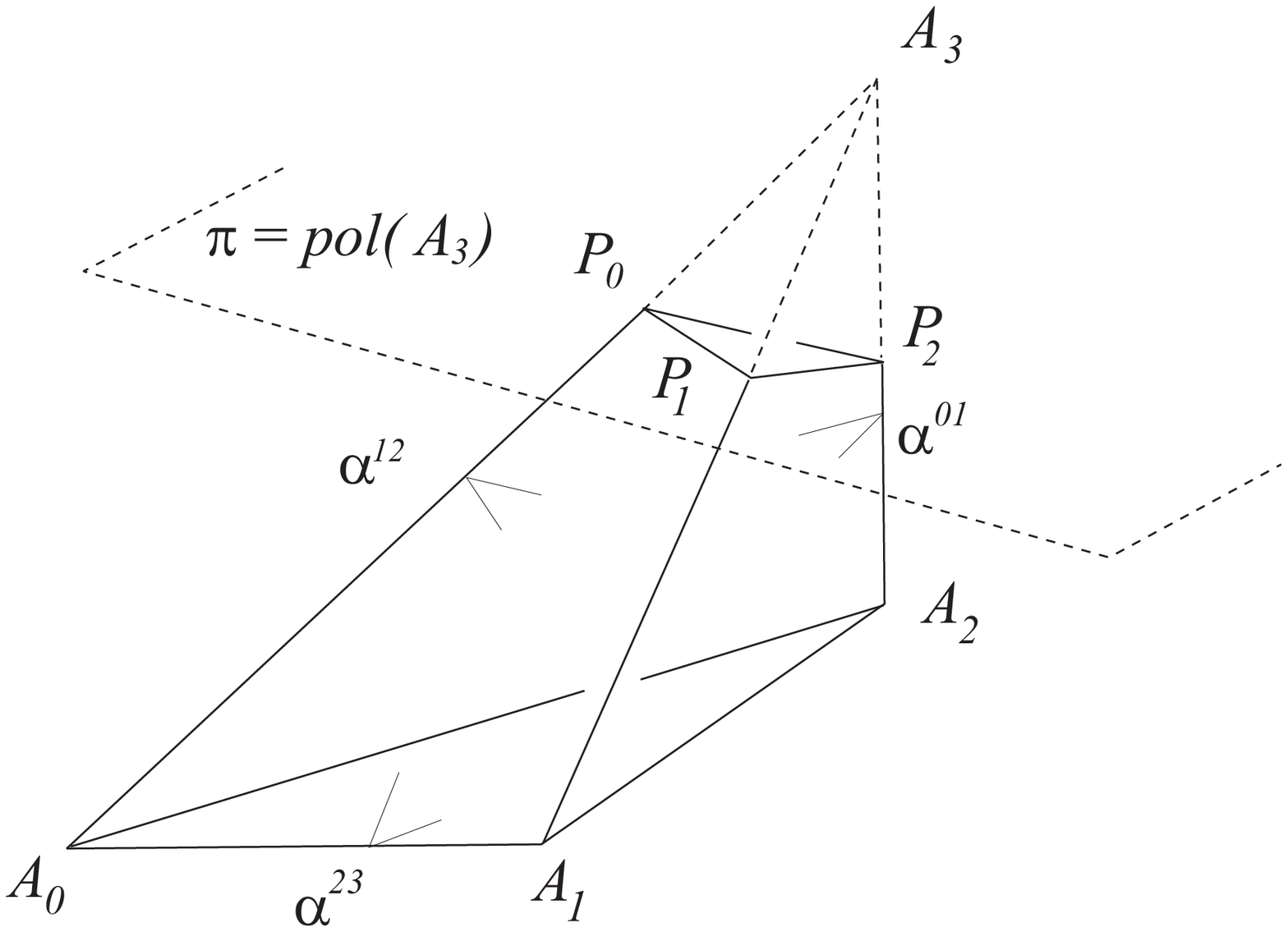} \includegraphics[width=6.5cm]{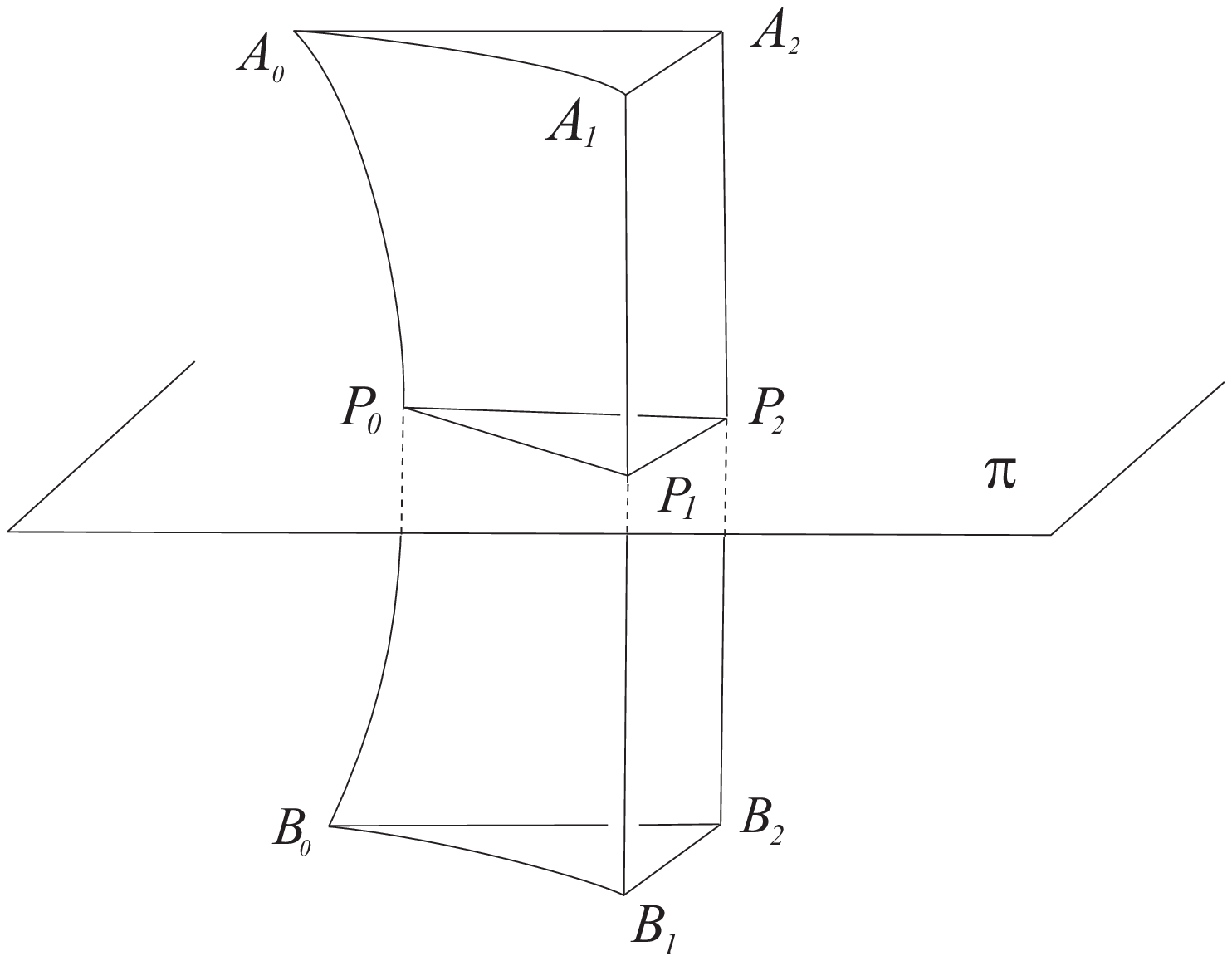}

a. \hspace{5cm} b.
\caption{a.~A $3$-dimensional complete orthoscheme of degree $d=1$ (simple frustum orthoscheme)
with outer vertex $A_3$. This orthoscheme is truncated by its polar plane $\pi=pol(A_3)$.~b.~Two congruent
adjacent simple frustum orthoschemes.}
\label{}
\end{figure}
We consider the images of $\mathcal{S}$ under reflections on its side facets.
The union of these $n$-dimensional orthoschames (having the common $\pi$ hyperplane) forms an infinite polyhedron denoted by $\mathcal{G}$.
$\mathcal{G}$ and its images under reflections on its ,,cover facets" fill hyperbolic
space $\mathbb{H}^n$ without overlap and generate $n$-dimensional tilings $\mathcal{T}$.

{\it The constant $k =\sqrt{\frac{-1}{K}}$ is the natural length unit in
$\mathbb{H}^n$. $K$ will be the constant negative sectional curvature. In the following we assume that $k=1$.}
\subsection{Volumes of the $n$-dimensional \\ Coxeter orthoschemes}

\begin{enumerate}
\item $2$-dimensional hyperbolic space $\mathbb{H}^2$

In the hyperbolic plane a simple frustum orthoscheme is a Lambert quadrilateral with exactly three right angles and its fourth angle is acute
$\frac{\pi}{q}$ ($q \ge 3$) (see Fig.~1). In our case the Lambert quadrilateral has a vertex at the infinity i.e. the angle at this vertex is $0$.
Its area can be determined by the well-known defect formula of hyperbolic triangles (see \cite{GHA}):
\begin{equation}
Vol_2(\mathcal{S})=\frac{\pi}{2}. \tag{3.1}
\end{equation}

\item $3$-dimensional hyperbolic space $\HYP$:

{Our polyhedron $A_0A_1A_2P_0P_1P_2$ is a simple frustum orthoscheme with
outer vertex $A_3$ (see Fig.~1) whose volume can be calculated by the following theorem of R.~Kellerhals
\cite{K89}:}
\begin{theorem} The volume of a three-dimensional hyperbolic
complete ortho\-scheme (except Lambert cube cases) $\mathcal{S}$
is expressed with the essential angles $\alpha_{01},\alpha_{12},\alpha_{23}, \ (0 \le \alpha_{ij} \le \frac{\pi}{2})$
(Fig.~1) in the following form:

\begin{align}
&Vol_3(\mathcal{S})=\frac{1}{4} \{ \mathcal{L}(\alpha_{01}+\theta)-
\mathcal{L}(\alpha_{01}-\theta)+\mathcal{L}(\frac{\pi}{2}+\alpha_{12}-\theta)+ \notag \\
&+\mathcal{L}(\frac{\pi}{2}-\alpha_{12}-\theta)+\mathcal{L}(\alpha_{23}+\theta)-
\mathcal{L}(\alpha_{23}-\theta)+2\mathcal{L}(\frac{\pi}{2}-\theta) \}, \tag{3.2}
\end{align}
where $\theta \in [0,\frac{\pi}{2})$ is defined by the following formula:
$$
\tan(\theta)=\frac{\sqrt{ \cos^2{\alpha_{12}}-\sin^2{\alpha_{01}} \sin^2{\alpha_{23}
}}} {\cos{\alpha_{01}}\cos{\alpha_{23}}}
$$
and where $\mathcal{L}(x):=-\int\limits_0^x \log \vert {2\sin{t}} \vert dt$ \ denotes the
Lobachevsky function.
\end{theorem}
For our prism tilings $\mathcal{T}_{pqr}$ we have:~
$\alpha_{01}=\frac{\pi}{p}, \ \ \alpha_{12}=\frac{\pi}{q}, \ \
\alpha_{23}=\frac{\pi}{r}$ .
\end{enumerate}
\subsection{On hyperballs}
The equidistant surface (or hypersphere) is a quadratic surface that lies at a constant distance
from a plane in both halfspaces. The infinite body of the hypersphere is called a hyperball.
The $n$-dimensional {\it half-hypersphere } $(n=2,3)$ with distance $h$ to a hyperplane $\pi$
is denoted by $\mathcal{H}_n^h$.
The volume of a bounded hyperball piece $\mathcal{H}_n^h(\mathcal{A}_{n-1})$
bounded by an $(n-1)$-polytope $\mathcal{A}_{n-1} \subset \pi$, $\mathcal{H}_n^h$ and by
hyperplanes orthogonal to $\pi$ derived from the facets of $\mathcal{A}_{n-1}$ can be determined by the
formulas (3.3) and (3.4) that follow from the suitable extension of the classical method of {{J.~Bolyai}}:
\begin{equation}
Vol_2(\mathcal{H}_2^h(\mathcal{A}_1))=Vol_1(\mathcal{A}_{1}) \sinh{(h)}, \tag{3.3}
\end{equation}
\begin{equation}
Vol_3(\mathcal{H}_3^h(\mathcal{A}_2))=\frac{1}{4}Vol_2(\mathcal{A}_{2})\left[\sinh{(2h)}+
2 h \right], \tag{3.4}
\end{equation}
where the volume of the hyperbolic $(n-1)$-polytope $\mathcal{A}_{n-1}$ lying in the plane
$\pi$ is $Vol_{n-1}(\mathcal{A}_{n-1})$.
\subsection{On horoballs}
A horosphere in $\mathbb{H}^n$ ($n \ge 2)$ is a
hyperbolic $n$-sphere with infinite radius centered
at an ideal point on $\partial \mathbb{H}^n$. Equivalently, a horosphere is an $(n-1)$-surface orthogonal to
the set of parallel straight lines passing through a point of the absolute quadratic surface.
A horoball is a horosphere together with its interior.

We consider the usual Beltrami-Cayley-Klein ball model of $\mathbb{H}^n$
centered at $O(1,0,0,$ $\dots, 0)$ with a given vector basis
$\bold{e}_i \ (i=0,1,2,\dots, n)$ and set an
arbitrary point at infinity to lie at $T_0=(1,0,\dots, 0,1)$.
The equation of a horosphere with center
$T_0=(1,0,\dots,1)$ passing through point $S=(1,0,\dots,s)$ is derived from the equation of the
the absolute sphere $-x^0 x^0 +x^1 x^1+x^2 x^2+\dots + x^n x^n = 0$, and the plane $x^0-x^n=0$ tangent to the absolute sphere at $T_0$.
The general equation of the horosphere is in projective coordinates ($s \neq \pm1$):
\begin{align}
(s-1)\left(-x^0 x^0 +\sum_{i=1}^n (x^i)^2\right)-(1+s){(x^0-x^n)}^2 & =0, \tag{3.5}
\end{align}
and in cartesian coordinates setting $h_i=\frac{x^i}{x^0}$ it becomes
\begin{equation}
\label{eqn:horosphere1}
\frac{2 \left(\sum_{i=1}^{n-1} h_i^2 \right)}{1-s}+\frac{4 \left(h_n-\frac{s+1}{2}\right)^2}{(1-s)^2}=1. \tag{3.6}
\end{equation}

In $n$-dimensional hyperbolic space any two horoballs are congruent in the classical sense.
However, it is often useful to distinguish between certain horoballs of a packing.
We use the notion of horoball type with respect to the packing as introduced in \cite{Sz12-2}.

In order to compute volumes of horoball pieces, we use J\'anos Bolyai's classical formulas from the mid 19-th century:
\begin{enumerate}
\item
The hyperbolic length $L(x)$ of a horospheric arc that belongs to a chord segment of length $x$ is
\begin{equation}
\label{eq:horo_dist}
L(x)=2 \sinh{\left(\frac{x}{2}\right)}. \tag{3.7}
\end{equation}
\item The intrinsic geometry of a horosphere is Euclidean,
so the $(n-1)$-dimensional volume $\mathcal{A}$ of a polyhedron $A$ on the
surface of the horosphere can be calculated as in $\mathbb{E}^{n-1}$.
The volume of the horoball piece $\mathcal{H}(A)$ determined by $A$ and
the aggregate of axes
drawn from $A$ to the center of the horoball is
\begin{equation}
\label{eq:bolyai}
Vol(\mathcal{H}(A)) = \frac{1}{n-1}\mathcal{A}. \tag{3.8}
\end{equation}
\end{enumerate}
\section{Hyp-hor packings in hyperbolic plane}
We consider the previously described $2$-dimensional Coxeter tilings given by the Coxeter symbol $[\infty]$ (see Fig.~1),
which are denoted by $\mathcal{T}_a$. The fundamental domain of $\mathcal{T}_a$ is a Lambert quadrilateral $A_0A_1P_1P_0$
(see Fig.~3) that is denoted by $\mathcal{F}_a$. It is derived by the truncation of the orthoscheme $A_0A_1A_2$
by the polar line $\pi$ of vertex $A_2$ where the initial principal vertex of the orthoschemes $A_0$ is lying on the absolute quadric $Q$
and its other principal vertex $A_2$ is an outer point of the model.

Its images under reflections on its sides fill hyperbolic
plane $\mathbb{H}^2$ without overlap. The tilings $\mathcal{T}_a$ contain a free parameter $0 < a < 1,~ ~ a \in \bR$. The
polar straight line of $A_2$ is $\pi$ and $\pi \cap A_0A_2=P_0$, $\pi \cap A_1A_2=P_1$.

We consider the usual Beltrami-Cayley-Klein ball model of $\mathbb{H}^2$
centered at $O(1,0,0)$ with a given vector basis
$\bold{e}_i \ (i=0,1,2)$ and set the above Lambert quadrilateral $A_0A_1P_1P_0$ in this coordinate system with coordinates
$$
A_0(1,0,1);~A_1(1,0,0);~P_1(1,a,0);~P_0(1,a,1-a^2); ~(0 < a < 1).
$$
The polar line $\Bu_2(1,-\frac{1}{a},0)^T$ of the outer vertex $A_2(1,\frac{1}{a},0)$ is $\pi$ which contains the points
$P_0$ and $P_1$ (see Fig.~3).

We construct hyp-hor packings to $\mathcal{T}_a$ tilings therefore the hyper- and horocycles have to satisfy the following requirements:
\begin{enumerate}
\item The centre of the horocycle can only be the vertex $A_0$ and the
corresponding horocycle $\mathfrak{H}_a(y_1)$ has not common points with inner of segments $A_1P_1$ and $P_0P_1$.
These horocycle types depend on parameter $a$ of the considered tiling $\mathcal{T}_a$ and passing through
the point $T(1,0,y_1)$ $(0< a,~y_1 <1)$ (see Fig.~4).
\item We can choose the base straight line of the hypercycle $\mathcal{H}_a(y_2)$ between the lines $P_0P_1$ and $A_1P_1$, the role of these
lines is symmetrical regarding the packings. We consider the $A_1P_1$ line as base line of hypercycles to construct hyp-hor packings.
Furthermore, $\mathcal{H}_a(y_2)$ has not common points with inner of segments $A_0P_0$. These hypercycle types
depend on the parameter parameter $a$ of the considered tiling $\mathcal{T}_a$ and passing through
the points $T(1,0,y_2)$ $(0< a,~y_2 <1)$ (see Fig.~4).
\item $ card \{ int (\mathfrak{H}_a(y_1))\cap int(\mathcal{H}_a(y_2))\}=0.$
\end{enumerate}
{\it If the hyper-and horocycles hold the above requirements then we obtain hyp-hor packings $\mathcal{T}_a$ in the hyperbolic plane
derived by the structure of the considered Coxeter simplex tilings.}
\begin{figure}[ht]
\centering
\includegraphics[width=9cm]{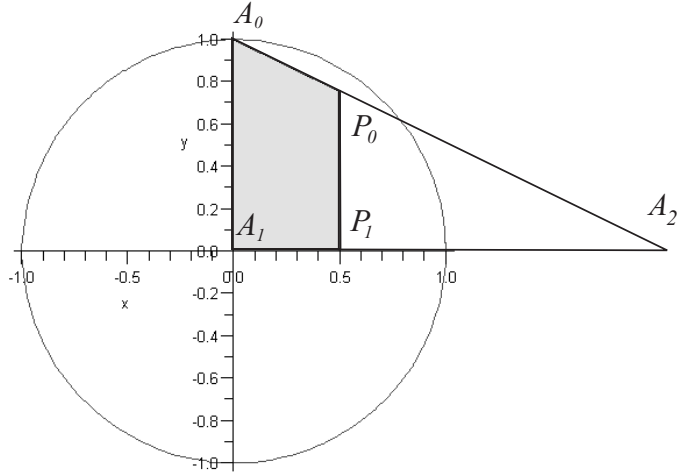}
\caption{The fundamental domain of $\mathcal{T}_a$ is a Lambert quadrilateral $A_0A_1P_0P_1$, at present $a=0.5$.}
\label{}
\end{figure}
\begin{definition}
The density of the above hyp-hor packings $\mathcal{P}_a(y_1,y_2)$ is
\begin{equation}
\delta(\mathcal{P}_a(y_1,y_2))=\frac{Vol(\mathcal{F}_a \cap (\mathfrak{H}_a(y_1) \cup \mathcal{H}_a(y_2))}{Vol(\mathcal{F}_a)}. \notag
\end{equation}
\end{definition}
It is well known that a packing is locally optimal (i.e. its density is locally maximal), then it is locally stable i.e. each ball is fixed by
the other ones so that no ball of packing
can be moved alone without overlapping another ball of the given ball packing or by other requirements of the corresponding tiling. Therefore, we can assume that the horocycle $\mathfrak{H}_a(y_1)$
and the hypercycle $\mathcal{H}_a(y_2)$ touch each other at the point $T(1,0,y)$ where $(0< y <1)$.
The possible values of $y=y_1=y_2$ may depend on parameter $a$ (see Fig.~4).
\begin{figure}[ht]
\centering
\includegraphics[width=6cm]{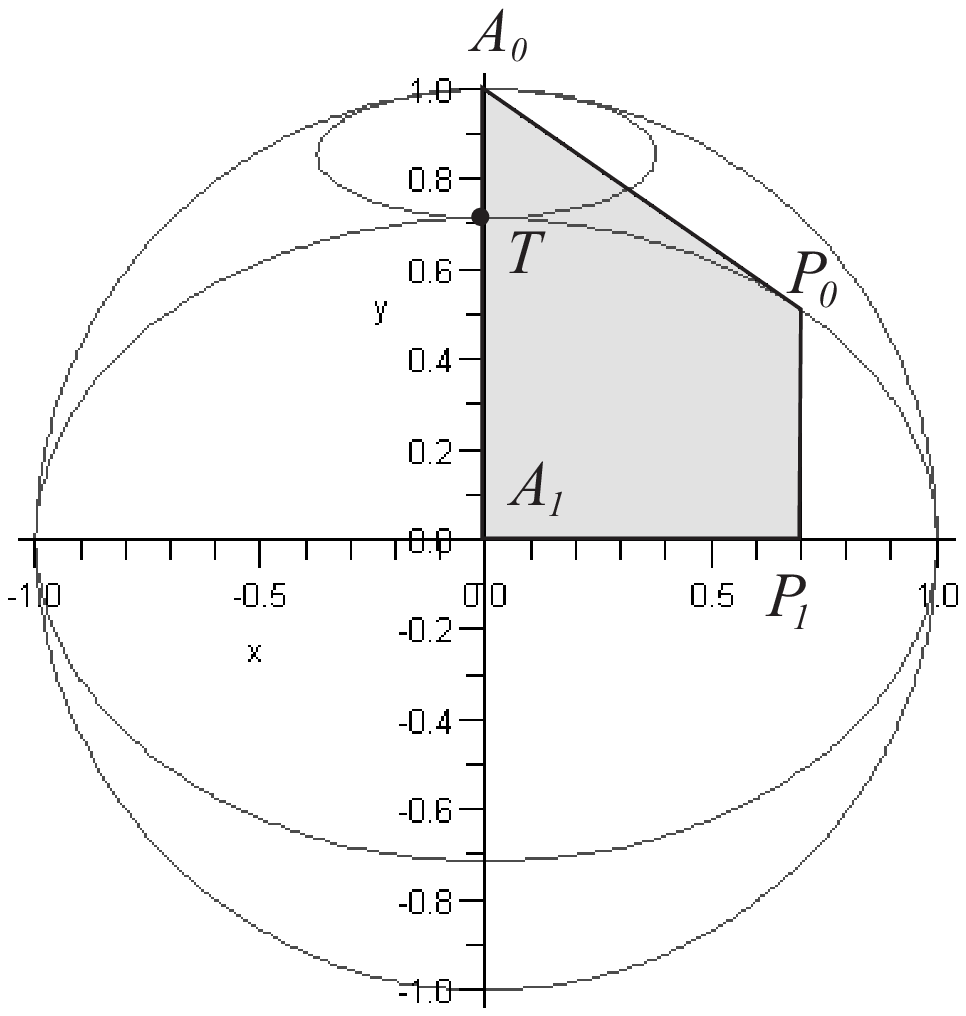} ~ ~ \includegraphics[width=6cm]{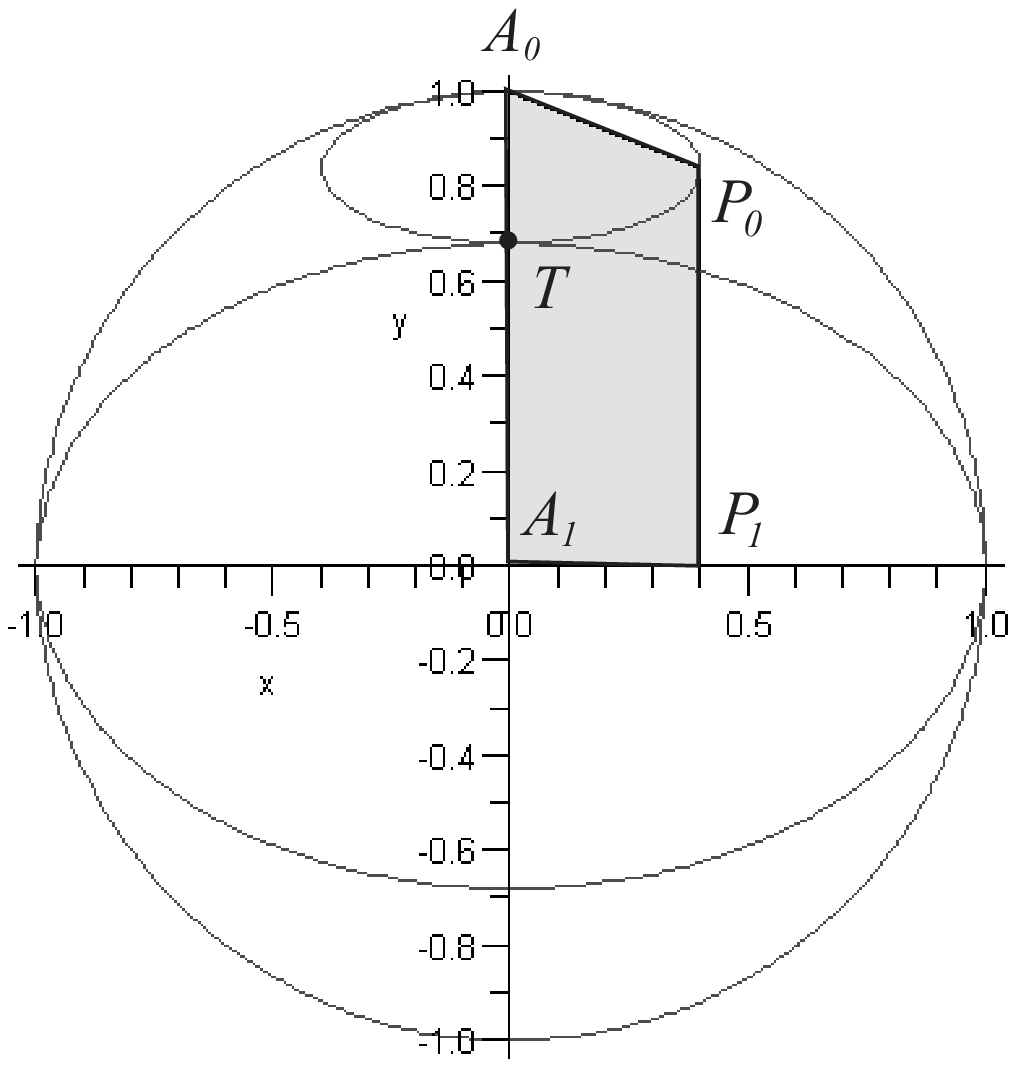}

a. \hspace{6cm} b.
\caption{a. The hyp-hor packing $\mathcal{P}_a^1$ of type 1, at present $a=0.7$. ~
b. The hyp-hor packing $\mathcal{P}_a^2$ of type 2, at present $a=0.4$.}
\label{}
\end{figure}
\subsection{Main types of hyp-hor packings}
We distinguish two main types of hyp-hor packings:
\begin{enumerate}
\item The hypercycle $\mathcal{H}_a(y)$ contains the point $P_0$ and the horocycle $\mathfrak{H}_a(y)$ touches it at $T(1,0,y)$. These
configuration can be realized for all possible parameters $0<a<1$. These packings are denoted by $\mathcal{P}_a^1(y)$
(see Fig~4.a).
\item The horocycle $\mathfrak{H}_a(y)$ passes through the point $P_0$ and the hypercycle $\mathcal{H}_a(y)$ touches it at $T(1,0,y)$.
These configurations exist if $0<a\le\frac{1}{\sqrt{2}}$. If $a=\frac{1}{\sqrt{2}}$ then $y=0$ i.e. $\mathfrak{H}_a(0)$ touches the line
$A_1P_1$ at $A_1$, (the height of the hypersphere is $0$). If $a > \frac{1}{\sqrt{2}}$ then these configurations do not satisfy the requirements
of the hyp-hor packings. These packings are denoted by $\mathcal{P}_a^2(y)$ (see  Fig.~4.b).
\end{enumerate}
\subsubsection{The densities of packings $\mathcal{P}_a^1(y)$}
In this case the coordinates of touching point $T(1,0,y)$ can be easily expressed as the function of parameter $a$:
$y=\sqrt{1-a^2}$.
We obtain by the formulas (3.1), (3.3), (3.6), (3.7), (3.8) and by Definition 4.1 that the density of the packings $\mathcal{P}_a^1(y)$ of
type 1 can be calculated by the following formula:
\begin{equation}
\begin{gathered}
\delta(\mathcal{P}_a^1(y))=\delta(\mathcal{P}^1(a))=\\
\left( 4\,\sinh \left( \frac{1}{2} \,{{\rm arccosh}} \left( \frac{1}{2} \,{\frac {2-2\,
\sqrt {1-{a}^{2}}-{a}^{2}+2\,{a}^{4}}{{a}^{4}}} \right)  \right) +2\,
\sqrt {1-{a}^{2}} \right)/{\pi }  \tag{4.1}
\end{gathered}
\end{equation}
where $0 < a < 1$.
\begin{lemma} Analysing the above density formula we obtain that
$$\lim_{a\rightarrow 0}{\Big[ \delta(\mathcal{P}^1(a))\Big]}=\frac{3}{\pi}, ~ \lim_{a\rightarrow 1}{\Big[\delta(\mathcal{P}^1(a))\Big]}=
\frac{2}{\pi}$$ and $\frac{2}{\pi} < \delta(\mathcal{P}^1(a))< \frac{3}{\pi}$ for parameters $0 < a < 1$ (see Fig.~5.a).
\end{lemma}
\begin{corollary}
In the hyperbolic plane $\mathbb{H}^2$ the universal upper bound density of ball packings
can be arbitrarily accurate approximate with the densities $\delta(\mathcal{P}^1(a))$ of hyp-hor packings of type 1.
\end{corollary}
\subsubsection{The densities of packings $\mathcal{P}_a^2(y)$}
Similarly to the previous section the coordinates of touching point $T(1,0,y)$ can be expressed as the function of parameter $a$:
$y=1-2 a^2$.
We obtain by the formulas (3.1), (3.3), (3.6), (3.7), (3.8) and by Definition 4.1 that the density of the packings $\mathcal{P}_a^2(y)$ of
type 2 can be calculated by the following formula:
\begin{equation}
\begin{gathered}
\delta(\mathcal{P}_a^2(x))=\delta(\mathcal{P}^2(a))= \\
\left( 4\,\sinh \left( \frac{1}{2} \,{\rm arccosh} \left( -\frac{1}{2}\,{\frac {-3+2\,
{a}^{2}}{1-{a}^{2}}} \right)  \right) -{\frac {-1+2\,{a}^{2}}{\sqrt {1
-{a}^{2}}}} \right) /{\pi } \tag{4.2}
\end{gathered}
\end{equation}
where $0 < a < \frac{1}{\sqrt{2}}$.
\begin{lemma} Analysing the above density formula we obtain that
$$\lim_{a\rightarrow 0}{\Big[ \delta(\mathcal{P}^2(a))\Big]}=\frac{3}{\pi}, ~ \lim_{a\rightarrow
\frac{1}{\sqrt{2}}}{\Big[\delta(\mathcal{P}^2(a))\Big]}=
\frac{2\sqrt{2}}{\pi}$$ and $\frac{2\sqrt{2}}{\pi} < \delta(\mathcal{P}^1(a))< \frac{3}{\pi}$ for parameters $0 < a < \frac{1}{\sqrt{2}}$
(see Fig.~5.b).
\end{lemma}
\begin{corollary}
In the hyperbolic plane $\mathbb{H}^2$ the universal upper bound density of ball packings
can be arbitrarily accurate approximate with the densities $\delta(\mathcal{P}^2(a))$ of hyp-hor packings of type 2.
\end{corollary}
\begin{figure}[ht]
\centering
\includegraphics[width=6cm]{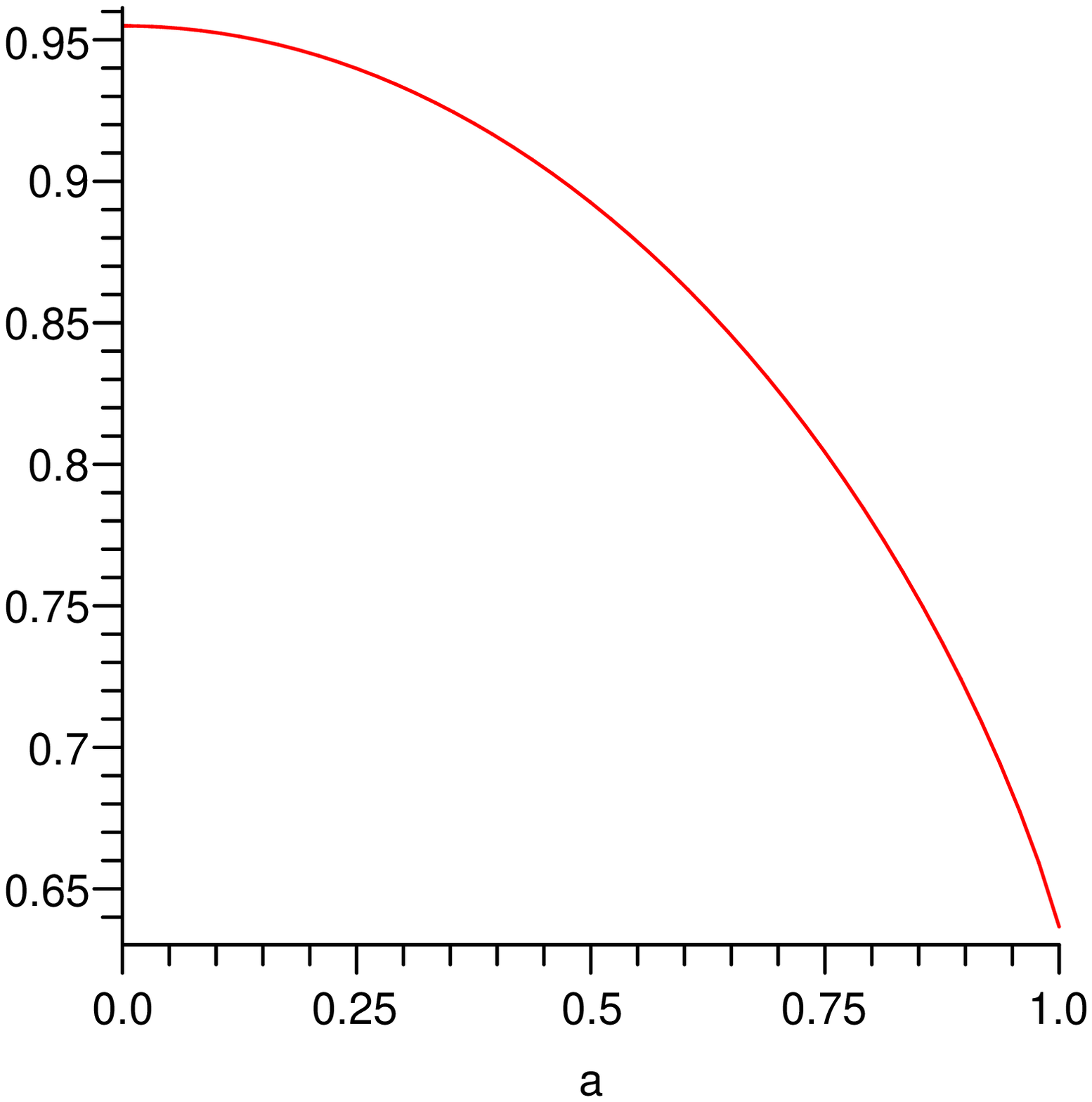} ~ ~ \includegraphics[width=6cm]{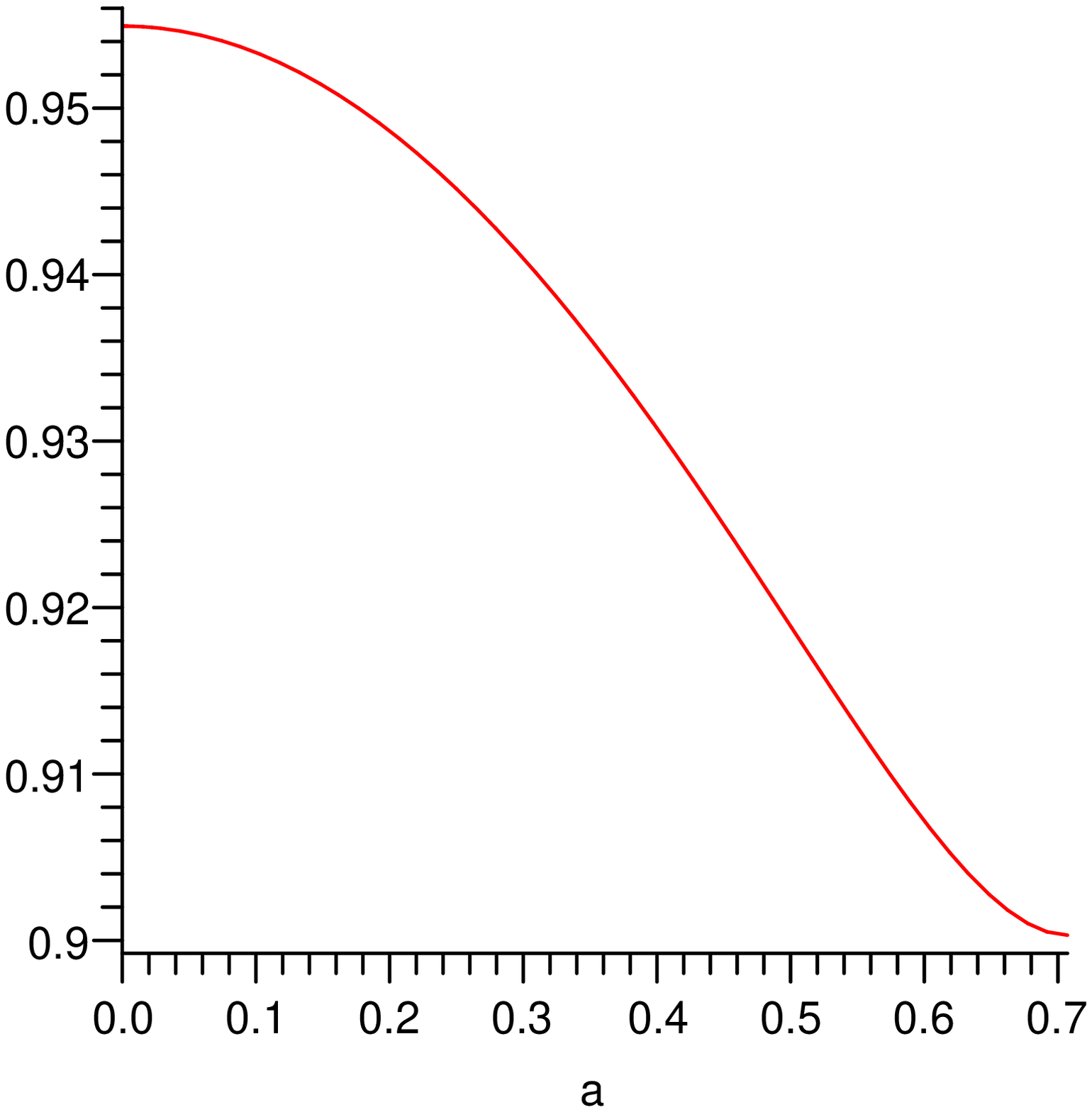}

a. \hspace{6cm} b.
\caption{a. The density function of hyp-hor packings $\mathcal{P}_a^1$ of type 1. ~
b. The density function of hyp-hor packings $\mathcal{P}_a^2$ of type 2.}
\label{}
\end{figure}
\subsection{The general cases}
\begin{enumerate}
\item First we consider the hyp-hor packings $\mathcal{P}_a(y)$ where the configurations are "between the two main cases": i.e. the inequalities
$0<a \le \frac{1}{\sqrt{2}}$ and $1-2a^2 < y \le \sqrt{1-a^2}$ hold.
\item We get the second case if the inequalities $\frac{1}{\sqrt{2}}< a < 1$ and  $0 \le y < \sqrt{1-a^2}$ hold.
\end{enumerate}
In both cases the densities of packings $\mathcal{P}_a(y)$ are denoted by
$\delta(\mathcal{P}_a(y))$ which can be determined by
the formulas (3.1), (3.3), (3.6), (3.7), (3.8) and by Definition 4.1:
\begin{equation}
\begin{gathered}
\delta(\mathcal{P}_a(y))=\Big( 4\,\sinh \left( \frac{1}{2}\,{\rm arccosh} \left( -\frac{1}{2}\,{\frac {-1+2\,
y-2\,{a}^{2}+2\,{a}^{2}{y}^{2}-{y}^{2}}{{a}^{2} \left( 1-{y}^{2}
\right) }} \right)  \right)\cdot \\ \cdot \sqrt {1-{y}^{2}}+2\,ya \Big) {\frac {1
}{\pi \sqrt {1-{y}^{2}}}} \tag{4.3}
\end{gathered}
\end{equation}
Analyzing the above density function we obtain that the maximal densities can be attained at the "main cases" described in
subsections 4.1.1 and 4.1.2. Therefore we get the following
\begin{theorem}
In the hyperbolic plane $\mathcal{H}^2$ the densest packing configurations with horo- and hyperballs generated by simple frustum orthoschemes
with Schl\"afly symbol $[\infty]$ provide the $\mathcal{P}^i(y)$ $(i=1,2)$ packings (described in subsection 4.1.1 and 4.1.2) if their
parameter $a \rightarrow 0$. Their densities $\delta(\mathcal{P}^i(a))$ are arbitrarily accurate approximate the
universal upper bound density $\frac{3}{\pi}$ of ball packings of $\mathcal{H}^2$.
\end{theorem}
\begin{rmrk}
If $y=0$ and $\frac{1}{\sqrt{2}} \le a < 1$ then we obtain ball packings $\mathcal{P}_a(0)$ which contain purely horocycles (see Fig.~6.a).
Their densities can be computed also by the formula (4.3) and its graph is illustrated in Fig.~7.  The maximal density is $0.90032$ belonging to
$a=\frac{1}{\sqrt{2}}$.
\end{rmrk}
\begin{figure}[ht]
\centering
\includegraphics[width=6cm]{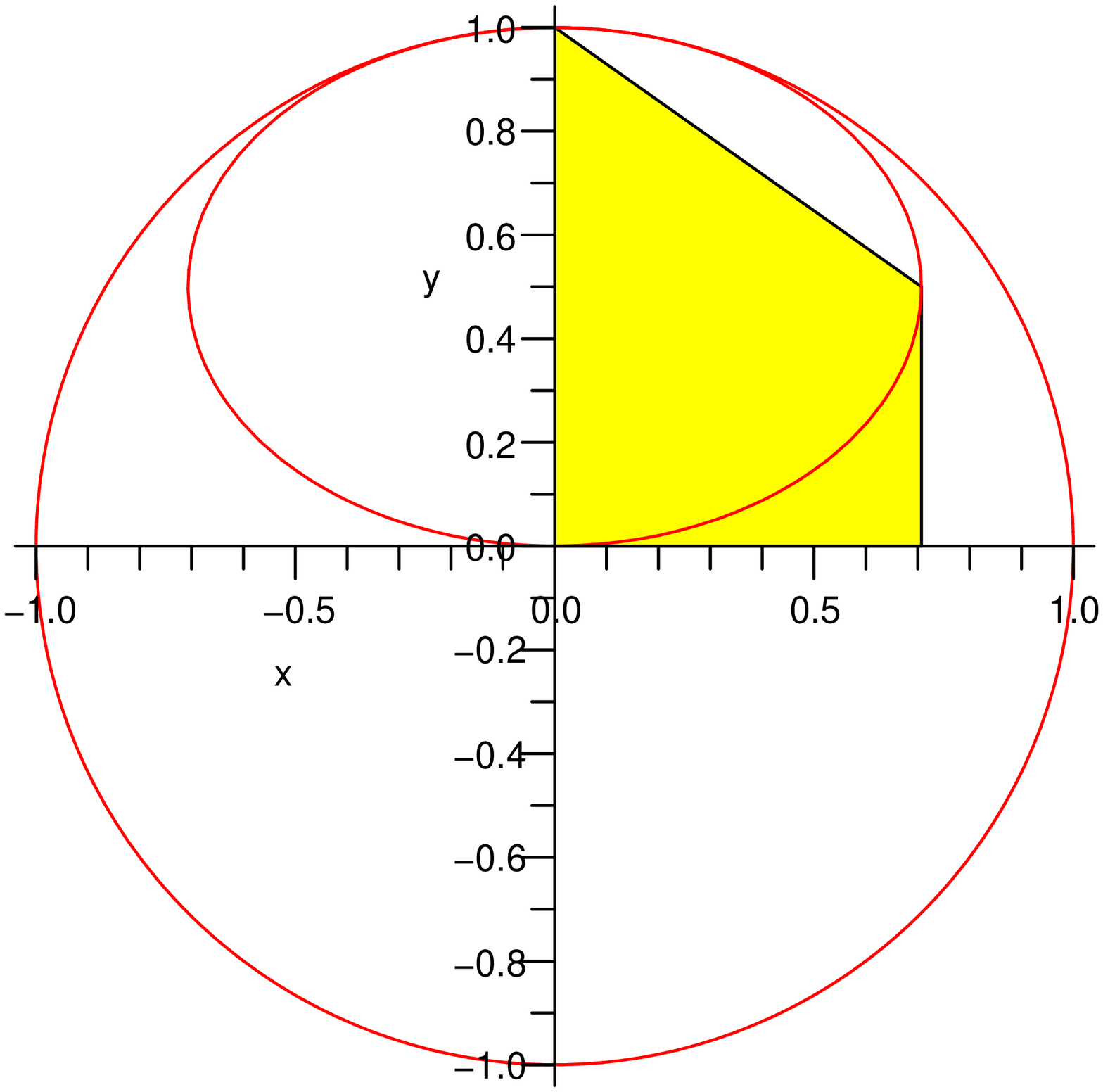} ~ ~ \includegraphics[width=6cm]{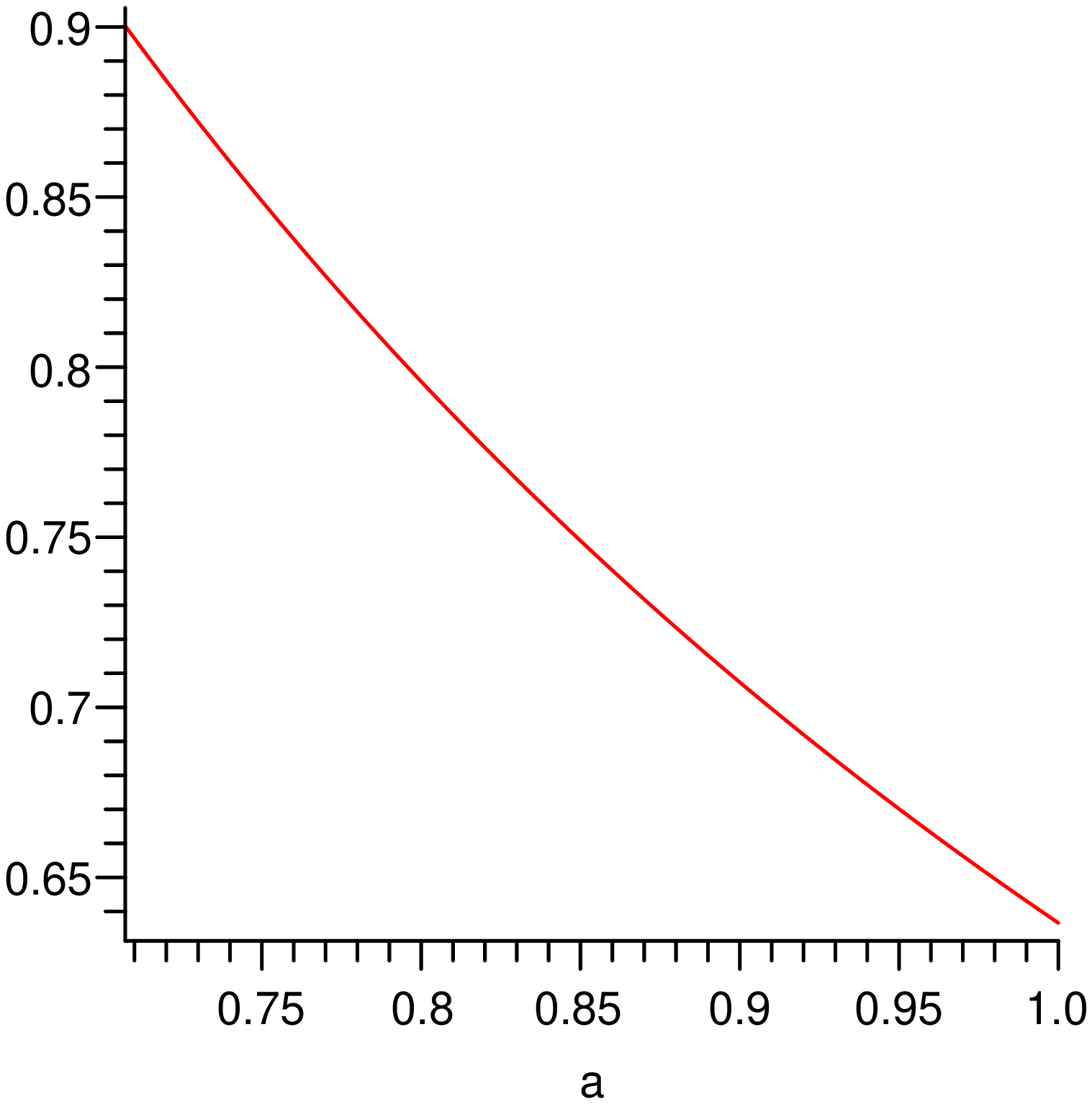}

a. \hspace{6cm} b.
\caption{a. The horocycle packing $\mathcal{P}_a(0)$ where $a=\frac{1}{\sqrt{2}}$. ~
b. The density function of horocycle packings $\mathcal{P}_a(0)$.}
\label{}
\end{figure}
\section{Hyp-hor packings in hyperbolic space $\HYP$}
In the 3-dimensional hyperbolic space there are $3$ infinite series of the simple frustum Coxeter orthoschemes with vertex at the infinity
listed in Fig.~1 and characterized in Sections 2-3.
The considered tilings with Schl\"afli symbol $[p,q,r]$, $(q,r)=(3,6), ~ (4,4), ~ (6,3)$ are denoted by $\mathcal{T}_p^{(i,j)}$ where
$p$ is an integer parameter and $p\ge 7$ if $(q,r)=(3,6)$, $p\ge 5$ if $(q,r)=(4,4)$, $p\ge 4$ if $(q,r)=(6,3)$.
The fundamental domain of $\mathcal{T}_p^{(q,r)}$ is a simple frustum orthoscheme $A_0A_1A_2P_0P_1P_2$ (see Fig.~2)
that is denoted by $\mathcal{F}_p^{(q,r)}$. It is derived by the truncation of the orthoscheme $A_0A_1A_2A_3$
by the polar hyperplane $\pi$ of vertex $A_3$ where the initial principal vertex of the orthoschemes $A_0$ is lying on the absolute quadric $Q$
and its other principal vertex $A_3$ is an outer point of the model.

Its images under reflections on its faces fill hyperbolic
space $\mathbb{H}^3$ without overlap. The polar plane of $A_3$ is $\pi$ and $\pi \cap A_iA_3=P_i$, $(i=0,1,2)$ (see Fig.~2).

We consider the usual Beltrami-Cayley-Klein ball model of $\mathbb{H}^3$
centered at $O(1,0,0,0)$ with a given vector basis
$\bold{e}_i \ (i=0,1,2,3)$ and set the above simple frustum orthoscheme $\mathcal{F}_p^{(q,r)}$ in this usual coordinate system (see Fig.~7.a-b).
\begin{equation}
\begin{gathered}
P_0(1,0,0,0);~P_1(1,0,y,0);~P_2(1,x,y,0)\ \ \text{where} \\
x=\sqrt{\tan^2{\frac{\pi}{q}}\frac{\cos^2{\frac{\pi}{q}}-\sin^2{\frac{\pi}{p}}}{\cos^2{\frac{\pi}{p}}\cos^2{\frac{\pi}{q}}}}; ~ ~
y=\sqrt{\frac{\cos^2{\frac{\pi}{p}}-\sin^2{\frac{\pi}{q}}}{\cos^2{\frac{\pi}{p}}}},\\
\text{and} ~ ~ A_0(1,0,0,1);~A_1(1,0,y,z_1);~A_2(1,x,y,z_2),~ \ \ \text{where} \\
\end{gathered} \tag{5.1}
\end{equation}
the $4{th}$ coordinates of the points $A_1$ and $A_2$ can be derived by the following procedure described in
general for $n$-dimensional hyperbolic space $\HYN$:
\begin{enumerate}
\item The points $P_k[{\mathbf{p}}_k]$ and $A_k[{\mathbf{a}}_k]$ $(k=1,2)$ are proper points of hyperbolic $n$-space and
$P_k$ lies on the polar hyperplane $pol(A_n)[\mbox{\boldmath$a$}^n]$ of the outer point $A_n$ thus
\begin{equation}
\begin{gathered}
\mathbf{p}_k \sim c \cdot \mathbf{a}_n+\mathbf{a}_k \in \mbox{\boldmath$a$}^n \Leftrightarrow
c \cdot \mathbf{a}_n \mbox{\boldmath$a$}^n+\mathbf{a}_k \mbox{\boldmath$a$}^n=0 \Leftrightarrow
c=-\frac{\mathbf{a}_k \mbox{\boldmath$a$}^n}{\mathbf{a}_n \mbox{\boldmath$a$}^n} \Leftrightarrow \\
\mathbf{p}_k \sim -\frac{\mathbf{a}_k \mbox{\boldmath$a$}^n}{\mathbf{a}_n \mbox{\boldmath$a$}^n}
\mathbf{a}_n+\mathbf{a}_k \sim \mathbf{a}_k (\mathbf{a}_n \mbox{\boldmath$a$}^n) - \mathbf{a}_n (\mathbf{a}_k \mbox{\boldmath$a$}^n)=
\mathbf{a}_k h_{nn}-\mathbf{a}_n h_{kn},
\end{gathered} \tag{5.2}
\end{equation}
where $h_{ij}$ is the inverse of the Coxeter-Schl\"afli matrix $c^{ij}$ (e.g. see (2.4) in $\mathbb{H}^3$) of the considered orthoscheme.
\item The hyperbolic distance $P_kA_k$ can be calculated by the following formula:
\[
\begin{gathered}
\cosh{P_kA_k}=\cosh{h}=\frac{- \langle {\mathbf{p}}_k, {\mathbf{a}}_k \rangle }
{\sqrt{\langle {\mathbf{p}}_k, {\mathbf{p}}_k \rangle \langle {\mathbf{a}}_k, {\mathbf{a}}_k \rangle}}= \\ =\frac{h_{kn}^2-h_{kk}h_{nn}}
{\sqrt{h_{kk}\langle \mathbf{p}_k, \mathbf{p}_k \rangle}} =
\sqrt{\frac{h_{kk}~h_{nn}-h_{kn}^2}
{h_{kk}~h_{nn}}}.
\end{gathered} \tag{5.3}
\]
\item The coordinates $z_k$ $(k=1,2)$ can be derived by the following equations (see Fig.~7.a-b):
\begin{equation}
\begin{gathered}
\cosh{P_1A_1}=\sqrt{\frac{h_{11}~h_{33}-h_{13}^2}{h_{11}~h_{33}}}=\frac{-1+x^2}{\sqrt{(-1+x^2)(-1+x^2+z_1^2)}},\\
\cosh{P_2A_2}=\sqrt{\frac{h_{22}~h_{33}-h_{23}^2}{h_{22}~h_{33}}}=\frac{-1+x^2+y^2}{\sqrt{(-1+x^2+y^2)(-1+x^2+y^2+z_2^2)}}.
\end{gathered} \tag{5.4}
\end{equation}
\end{enumerate}
For example for tiling $\mathcal{T}_7^{(3,6)}$ the above coordinates are the following (see Fig.~7.a-b):
{\small
\begin{equation}
\begin{gathered}
x=\frac{1}{18}\,{\frac {\sqrt {3}\sqrt {108 \,\cos^2 \left(\frac{\pi}{7}\,\right) -81}}{\cos\left( \frac{\pi}{7}
\right) }}\approx 0.27580, ~
y= \frac{1}{2}\,{\frac {\sqrt {3}\sqrt {-3+4\,  \cos^2 \left( \frac{\pi}{7}
\right) }}{\cos \left( \frac{\pi}{7}  \right) }}\approx 0.47770,\\
z_1= 12\,\sqrt {2}\sqrt {-566\,  \cos^2 \left( \frac{\pi}{7}  \right) +318+157\,\cos \left( \frac{\pi}{7}  \right) } \cos^2
\left( \frac{\pi}{7}  \right) \approx 0.92394 ,\\
z_2= 12\, \left( -1+2\, \cos^2 \left( \frac{\pi}{7}  \right) -
\cos \left( \frac{\pi}{7}  \right)  \right) \sqrt {18-31\,  \cos^2
\left( \frac{\pi}{7}  \right)  +8\,\cos \left( \frac{\pi}{7} \right) } \approx 0.69574.
\end{gathered} \notag
\end{equation}}
\begin{figure}[ht]
\centering
\includegraphics[width=6cm]{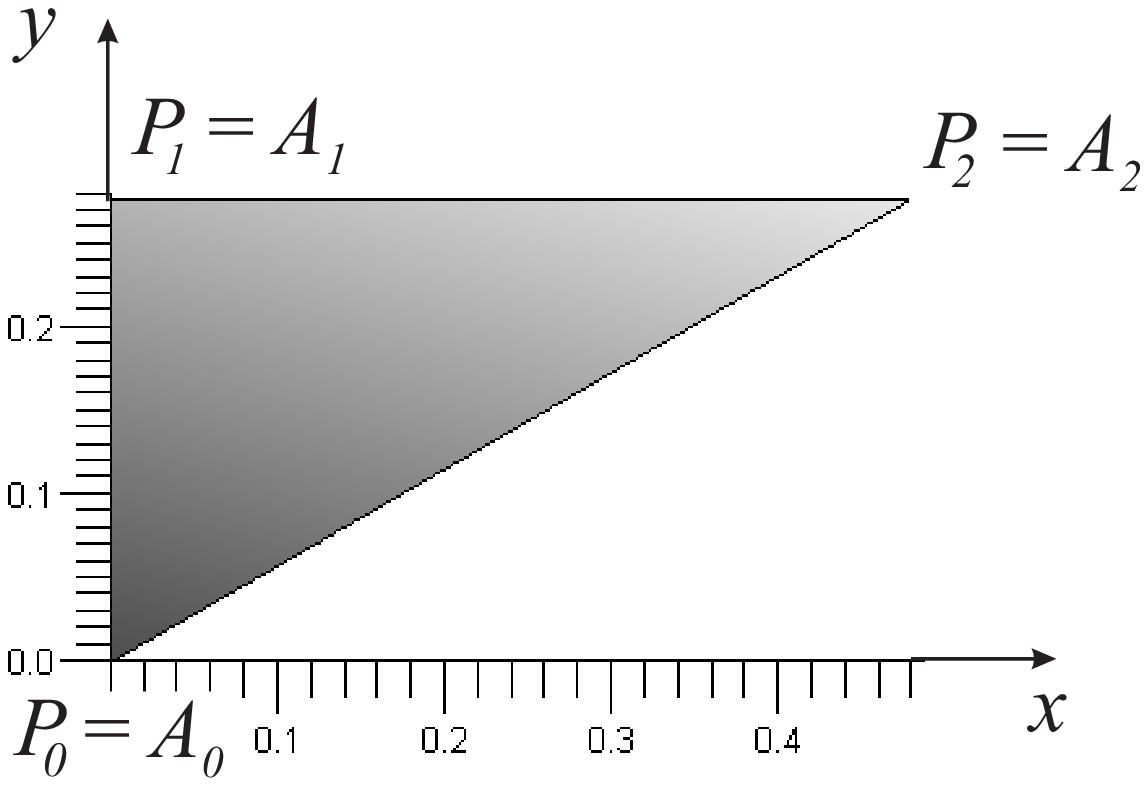} ~ ~ \includegraphics[width=6cm]{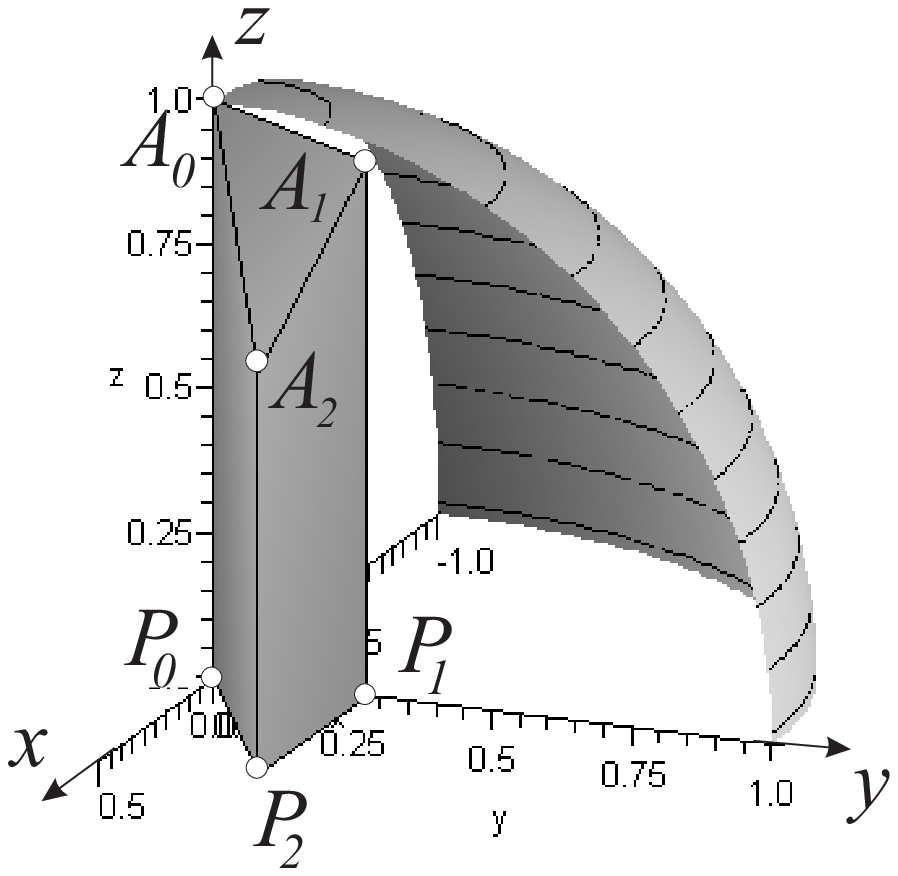}

a. \hspace{6cm} b.
\caption{The simple frustum orthoscheme with Schl\"afli symbol $[7,3,6]$ in Beltrami-Cayley-Klein model.}
\label{}
\end{figure}
We construct hyp-hor packings to $\mathcal{T}_p^{(q,r)}$ tilings therefore the hyper- and horocycles have to satisfy the following requirements:
\begin{enumerate}
\item The centre of the horoball can only be the vertex $A_0$ and the
corresponding horoball (horosphere) $\mathfrak{H}_p^{(q,r)}(\gamma_1(p))$ has not common points with inner of side faces $P_0P_1P_2$ and $P_1P_2A_1A_2$
of the simple frustum orthoschem $\mathcal{F}_p^{(q,r)}=P_0P_1P_2A_0A_1A_2$ (see Fig.~7.a-b).
These horoball types depend on the metric data of the considered tiling $\mathcal{T}_p^{(q,r)}$ and it is passing through 
the point $T_1(1,0,0,\gamma_1(p))$ where (see (5.3) and (5.4))
$0 \le \gamma_1(p) < 1.$
\item The base plane of the hyperball (hypersphere) $\mathcal{H}_p^{(q,r)}(\gamma_2(p))$ is $P_1P_2P_3$ plane.
Furthermore, the above hyperball has not common points with inner of side faces $A_0A_1A_2$. These hyperball types
depend on the metric data of the considered tiling $\mathcal{T}_p^{(q,r)}$ and it is passing through
the point $T_2(1,0,0,\gamma_2(p))$ where $0 \le \gamma_1(p) \le \tanh({\rm arcosh}(P_2A_2))<1$ (see (5.3) and (5.4))
(in some cases the points $T_1$ and $T_2$ can coincide i.e. the corresponding
horoball and hyperball touch each other).
\item $ card \{ int (\mathfrak{H}_p^{(q,r)}(\gamma_1(p)))\cap int(\mathcal{H}_p^{(q,r)}(\gamma_2(p)))\}=0.$
\end{enumerate}
{\it If the hyper-and horocycles hold the above requirements then we obtain hyp-hor packings $\mathcal{P}_p^{(q,r)}(\gamma_1(p),\gamma_2(p))$
in the hyperbolic plane derived by the structure of the considered Coxeter simplex tilings.}
\begin{definition}
The density of the above hyp-hor packings $\mathcal{P}_p^{(q,r)}(\gamma_1(p),\gamma_2(p))$ is
\begin{equation}
\delta(\mathcal{P}_p^{(q,r)}(\gamma_1(p),\gamma_2(p)))=\frac{Vol(\mathcal{F}_p^{(q,r)}
\cap (\mathfrak{H}_p^{(q,r)}(\gamma_1(p)) \cup \mathcal{H}_p^{(q,r)}(\gamma_2(p)))}{Vol(\mathcal{F}_p^{(q,r)})}. \notag
\end{equation}		
\end{definition}
\subsection{Hyp-hor packings to $\mathcal{T}_p^{(4,4)}$ and $\mathcal{T}_p^{(6,3)}$ tilings}
\subsubsection{On $\mathcal{T}_p^{(4,4)}$ tilings}
It is well known that a packing is locally optimal (i.e. its density is locally maximal), then it is locally stable i.e. each ball is fixed by
the other ones so that no ball of packing
can be moved alone without overlapping another ball of the given ball packing or by other requirements of the corresponding tiling.
Therefore, first we consider the largest possible horo- and hyperballs to the considered tiling.

The largest possible horoball centered at $A_0$ is passing through the vertex $A_1$ and the largest possible hyperball contains the vertex $A_2$.
We get by easy calculations, that these "maximal large balls" have common inner points for any permissible parameters $p$.
Thus, the optimal arrangement can be achieved if the horo- and hyperballs touch each other i.e. holds the $z(p)=\gamma_1(p)=\gamma_2(p)$ equation and
its density of the hyp-hor packing $\mathcal{P}_p^{(4,4)}(z(p))$ depends on parameter $z(p)$ where $0 \le z(p) \le \tanh({\rm arcosh}(P_2A_2))<1$
(see ((5.3) and (5.4)).
The volumes of $Vol(\mathcal{F}_p^{(4,4)} \cap \mathfrak{H}_p^{(4,4)})(z(p))$ and $Vol(\mathcal{F}_p^{(4,4)} \cap \mathcal{H}_p^{(4,4)})(z(p))$
can be calculated by the formulas (3.2), (3.4), (3.6), (3.7), (3.8), (5.3), (5.4) for any parameters $\bR \ni p > 4$ but only to parameters
$\bZ \ni p \ge 5$ belong tilings and so hyp-hor packings in $\HYP$.
The largest possible horoball is denoted by $\mathfrak{H}_p^{(4,4)}(z_{max}(p))$ passes through the point $T_1=T_2(1,0,0,z_{max}(p))$
which is the common point with the hypersphere $\mathcal{H}_p^{(4,4)}(z_{max}(p))$. We blow up this hypersphere
keeping the horoball $\mathfrak{H}_p^{(4,4)}(z_{max}(p))$ tangent to it upto this hypersphere touches the faces $A_1A_2A_3$ at vertex $A_2$.
During this process we can compute the densities of considered packings as the function of
$0 \le \zeta(p) \le {\rm arcosh}(P_2A_2))-{\rm artanh}(z_{max}(p))$ (see (5.4)) by the following equation
\begin{equation}
\begin{gathered}
\delta(\mathcal{P}_p^{(4,4)}(\zeta(p)))=\frac{Vol(\mathcal{F}_p^{(4,4)}
\cap (\mathfrak{H}_p^{(4,4)}(\zeta(p)) \cup \mathcal{H}_p^{(4,4)}(\zeta(p)))}{Vol(\mathcal{F}_p^{(4,4)})}, \\
{\rm artanh}(z(p))={\rm artanh}(z_{max}(p))+\zeta(p).
\notag
\end{gathered}
\end{equation}	
For example, if $p=5$ then $0  \le \zeta(5) \le \tanh({\rm arcosh}(P_2A_2))-z_{max}(5)  \approx 0.33419$ and the graph of
$\delta(\mathcal{P}_5^{(4,4)}(\zeta(5)))$ is described in Fig.~8.a. Analyzing the above density function we get that the
maximal density is achieved at the endpoint of the above interval with
density $0.81296$ (see Table 1).
\begin{figure}[ht]
\centering
\includegraphics[width=6cm]{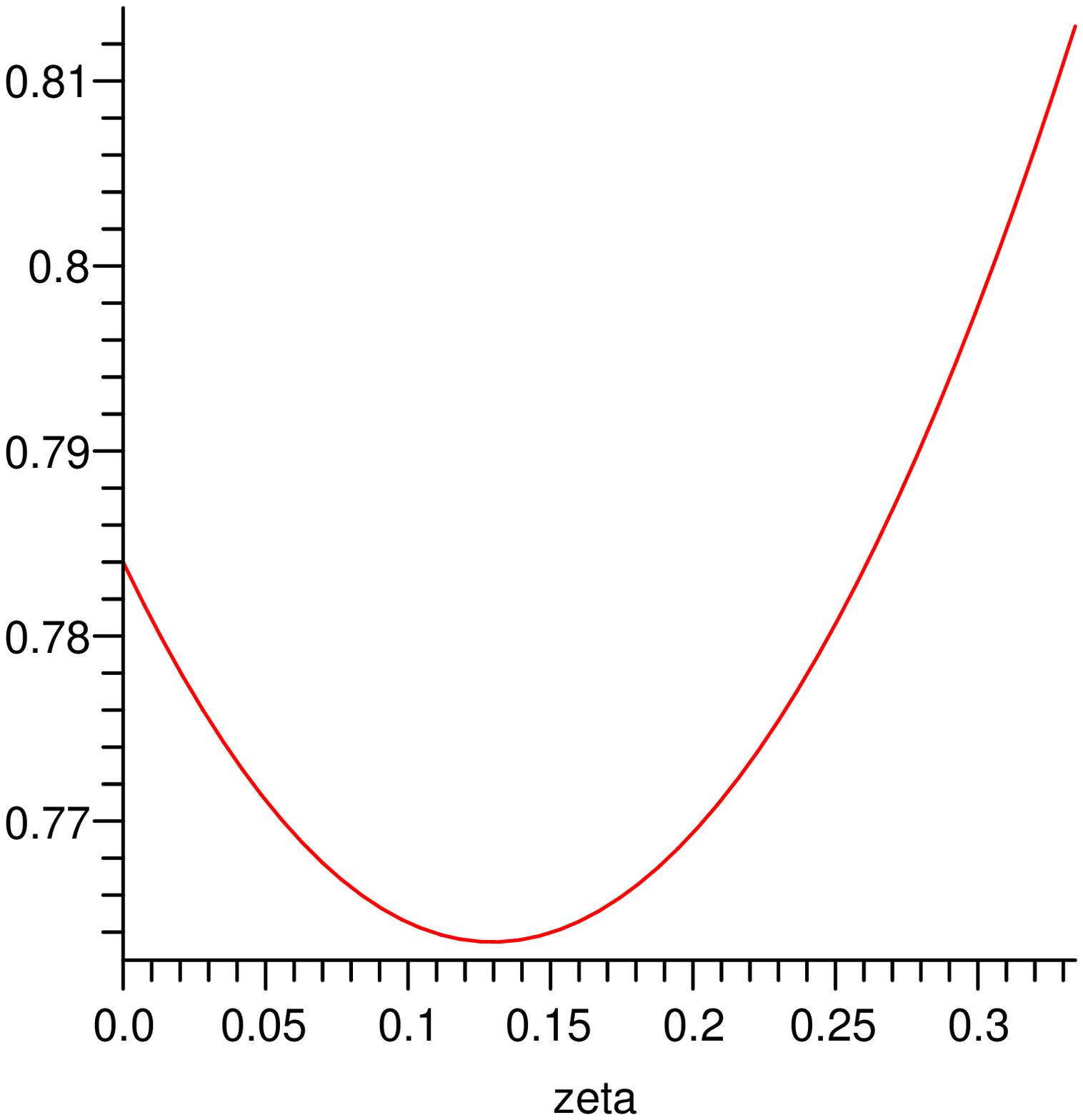} ~ ~ \includegraphics[width=6cm]{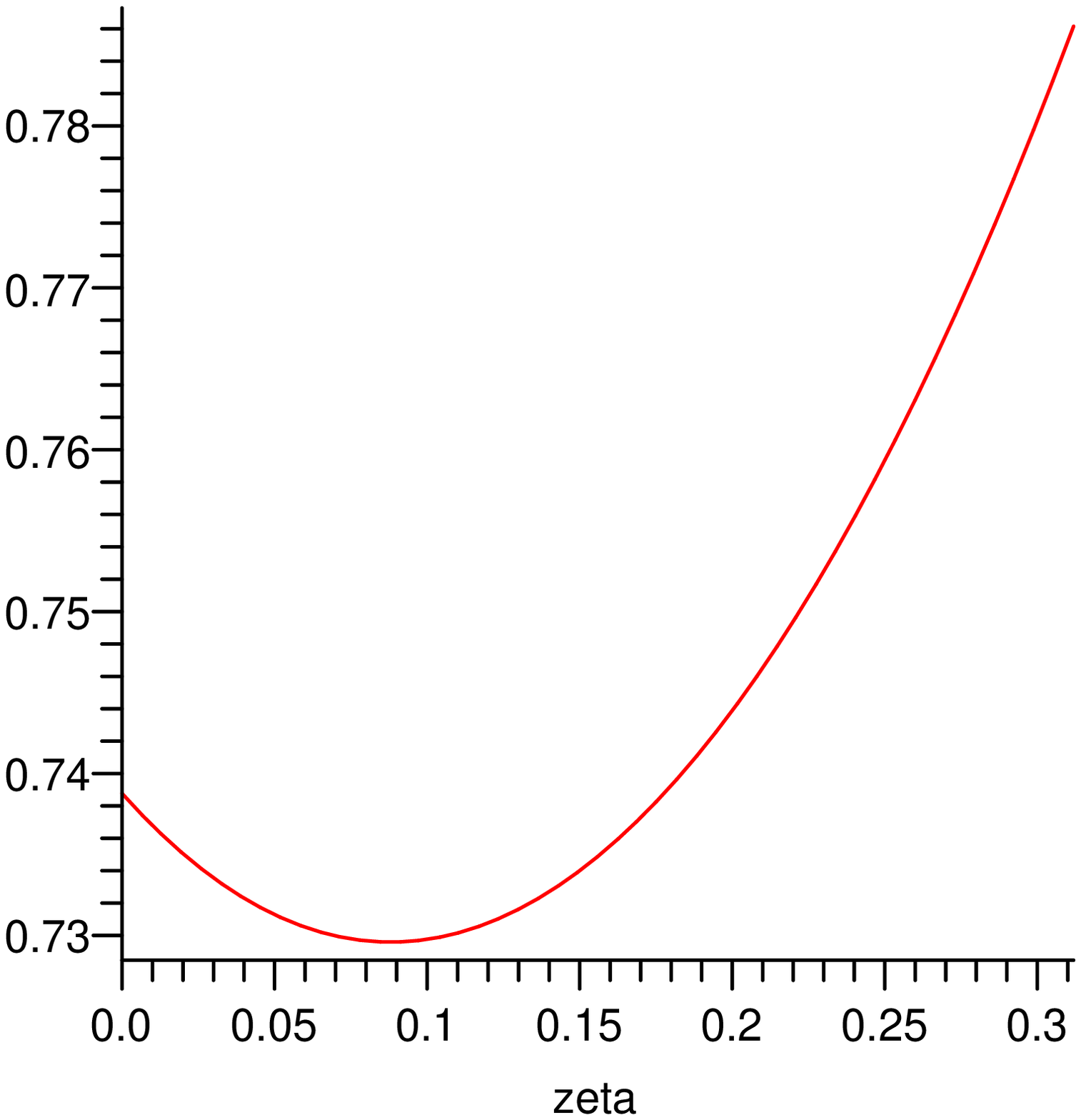}

a. \hspace{5cm} b.
\caption{The graphs of density functions $\delta(\mathcal{P}_5^{(4,4)}(\zeta(5)))$ and $\delta(\mathcal{P}_6^{(4,4)}(\zeta(6)))$.}
\label{}
\end{figure}
Similarly to the above $p=5$ parameter we can compute the optimal densities for all possible other parameters $p$.
The results for some parameters are summarized in Table 1.

{\it The densest hyp-hor configuration among the considered "realizable" packings belongs to packing $\mathcal{P}_{5}^{(4,4)}$
with density $\approx 0.812956$.}
\medbreak
{
\centerline{\vbox{
\halign{\strut\vrule~\hfil $#$ \hfil~\vrule
&\quad \hfil $#$ \hfil~\vrule
&\quad \hfil $#$ \hfil\quad\vrule
&\quad \hfil $#$ \hfil\quad\vrule
\cr
\noalign{\hrule}
\multispan4{\strut\vrule\hfill\bf Table 1, $q=4$, $r=4$  \hfill\vrule}%
\cr
\noalign{\hrule}
\noalign{\vskip2pt}
\noalign{\hrule}
p & Vol(\mathcal{F}_{p}^{(4,4)}) & {Vol(\mathcal{F}_p^{(4,4)}
\cap (\mathfrak{H}_p^{(4,4)} \cup \mathcal{H}_p^{(4,4)})})& \delta(\mathcal{P}_p^{(4,4)}) \cr
\noalign{\hrule}
5 & 0.34084197 & 0.27709010 &  0.81295769 \cr
\noalign{\hrule}
6 & 0.38165233 & 0.30003810 &  0.78615556 \cr
\noalign{\hrule}
7 & 0.40369221 & 0.30777518 &  0.76240058 \cr
\noalign{\hrule}}}}
}

\subsubsection{On $\mathcal{T}_p^{(6,3)}$ tilings}
The determination of the optimal hyp-hor packing configurations of packings $\mathcal{T}_p^{(6,3)}$ $(\bZ \ni p \ge 4)$ is similar to the above
tilings therefore here we only summarize the results in Table 2.

{\it The densest hyp-hor configuration among the considered packings belongs to packing $\mathcal{P}_{5}^{(4,4)}$
with density $\approx 0.81209$.}
\medbreak
{
\centerline{\vbox{
\halign{\strut\vrule~\hfil $#$ \hfil~\vrule
&\quad \hfil $#$ \hfil~\vrule
&\quad \hfil $#$ \hfil\quad\vrule
&\quad \hfil $#$ \hfil\quad\vrule
\cr
\noalign{\hrule}
\multispan4{\strut\vrule\hfill\bf Table 2, $q=6$, $r=3$  \hfill\vrule}%
\cr
\noalign{\hrule}
\noalign{\vskip2pt}
\noalign{\hrule}
p & Vol(\mathcal{F}_{p}^{(6,3)}) & {Vol(\mathcal{F}_p^{(6,3)}
\cap (\mathfrak{H}_p^{(6,3)} \cup \mathcal{H}_p^{(6,3)})})& \delta(\mathcal{P}_p^{(6,3)}) \cr
\noalign{\hrule}
4 & 0.31716925 & 0.25756985 &  0.81208961 \cr
\noalign{\hrule}
5 & 0.35991902 & 0.27187731 &  0.75538469 \cr
\noalign{\hrule}
6 & 0.38060310 & 0.27009741 &  0.70965634 \cr
\noalign{\hrule}}}}
}
\subsection{Hyp-hor packings to $\mathcal{T}_p^{(3,6)}$ tilings}
The investigations of these tilings are a little different from the above tilings.

First we consider the largest possible horo- and hyperballs to the considered tiling.

The largest possible horoball centered at $A_0$ is passing through the vertex $A_1$ and the largest possible hyperball contains the vertex $A_2$.
Contrary to the above $\mathcal{T}_p^{(4,4)}$ tilings here we get by easy calculations, that these
"maximal large balls" have not common points for any permissible parameters $p$.

The volumes of $Vol(\mathcal{F}_p^{(3,6)} \cap \mathfrak{H}_p^{(3,6)})(\gamma_1(p))$ and $Vol(\mathcal{F}_p^{(3,6)} \cap \mathcal{H}_p^{(3,6)})(\gamma_2(p))$
can be calculated by the formulas (3.2), (3.4), (3.6), (3.7), (3.8), (5.3), (5.4) for any parameters $\bR \ni p > 6$.
\begin{rmrk}
We note here, that only the parameters $ \bZ \ni p \ge 7$ provide
tilings and hyp-hor ball packing configurations in the hyperbolic space $\HYP$.
The other parameters provide locally optimal density as well but these hyp-hor
packing configurations can not be extended to the entirety of hyperbolic space $\mathbb{H}^3$.
\end{rmrk}
We can compute the densities for all possible parameters $p$ of hyp-hor packings. The results for some parameters are summarized in Table 1.

{\it The densest hyp-hor configuration among the considered "realizable" packings belongs to packing $\mathcal{P}_{7}^{(3,6)}$
with density $\approx 0.83267$.}
\medbreak
{
\centerline{\vbox{
\halign{\strut\vrule~\hfil $#$ \hfil~\vrule
&\quad \hfil $#$ \hfil~\vrule
&\quad \hfil $#$ \hfil\quad\vrule
&\quad \hfil $#$ \hfil\quad\vrule
\cr
\noalign{\hrule}
\multispan4{\strut\vrule\hfill\bf Table 3, $q=3$, $r=6$  \hfill\vrule}%
\cr
\noalign{\hrule}
\noalign{\vskip2pt}
\noalign{\hrule}
p & Vol(\mathcal{F}_{p}^{(3,6)}) & {Vol(\mathcal{F}_p^{(3,6)}
\cap (\mathfrak{H}_p^{(3,6)} \cup \mathcal{H}_p^{(3,6)})})& \delta(\mathcal{P}_p^{(3,6)}) \cr
\noalign{\hrule}
7 & 0.31781165 & 0.26463185 &  0.83266882 \cr
\noalign{\hrule}
8 & 0.34695830 & 0.27901923 &  0.80418664 \cr
\noalign{\hrule}
9 & 0.36482363 & 0.28351212 &  0.77712105 \cr
\noalign{\hrule}}}}
}%
Finally, we obtain by careful computations and investigations from the above method and results the following
\begin{figure}[ht]
\centering
\includegraphics[width=12cm]{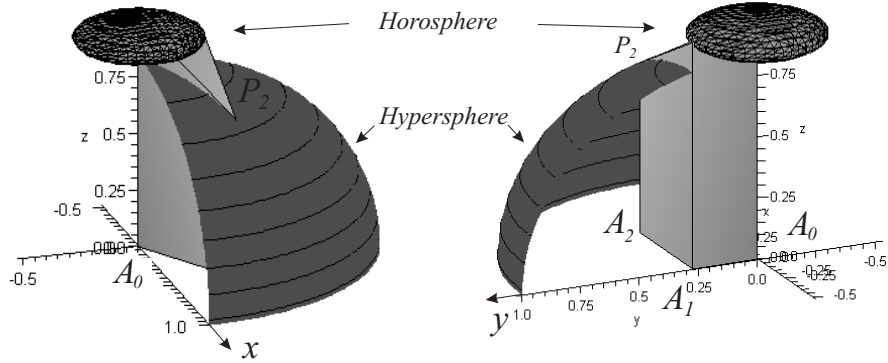}
\caption{The optimal hyp-hor packing configuration to Coxeter tiling $[7,3,6]$.}
\label{}
\end{figure}
\begin{theorem}
The $\mathcal{P}_7^{(3,6)}$ packing configuration (see Section 5.2) provides the maximal density $\approx 0.83267$
of hyp-hor packings $\mathcal{P}_p^{(q,r)})$ ($(q,r)=(4,4),$ $(6,3),(3,6)$ and
$p$ is suitable integer parameter (see Fig.~1)) which are derived by the Coxeter tilings generated by complete orthoschemes of degree 1
(simple frustum orthoschemes).
\end{theorem}
\subsubsection{On non-extendable hyp-hor packings $\mathcal{P}_p^{(3,6)}$ $(6 < p < 7, ~ p\in \bR)$}
The computation method described in the former sections is suitable to determine the densities of hyp-hor packings for
parameters $(6 < p < 7, ~ p\in \bR)$ as well. To any parameter belongs a simple frustum orthoscheme and therefore we can determine similarly to the
above cases the corresponding density of its optimal hyp-hor packing. But these packings can not be extended to the 3-dimensional
space. Analyzing these non-extendable packings for parameters $(6 < p < 7, ~ p\in \bR)$ we obtain the following
\begin{theorem}
The function $\delta(\mathcal{P}_p^{(3,6)})$, $(6 < p < 7, ~ p\in \bR)$
is attained its maximum for the parameter $p_{opt}$ which lies in the interval $[6.05,6.06]$
and the densities for parameters lying in this interval are larger that $\approx 0.85397$. That means that these
locally optimal hyp-hor configurations provide larger densities that the B\"or\"oczky-Florian density upper bound $(\approx 0.85328)$ for ball and
horoball packings (\cite{B--F64}).
\end{theorem}
\begin{rmrk}
{We note here, that the $5$-dimensional analogous periodic hyp-hor packing will be investigated in a forthcoming paper.}
\end{rmrk}
The question of finding the densest hyp-hor packing without any symmetry assumption in the $n$-dimensional hyperbolic space is open.
Similarly to it, the discussion of the densest horoball and hyperball packings in the $n$-dimensional hyperbolic space $n \ge 3$ with horoballs
of different types and congruent hyperballs has not been settled yet (see \cite{KSz}, \cite{KSz14}, \cite{Sz12}, \cite{Sz12-2}).

Moreover, optimal sphere packings in other homogeneous Thurston geometries represent
another huge class of open mathematical problems. For these non-Euclidean geometries
only very few results are known (e.g. \cite{Sz07-2}, \cite{Sz10}, \cite{Sz13-2}, \cite{Sz14-1}).
By the above these we can say that the revisited Kepler problem keep several interesting open questions.


\noindent
\footnotesize{Budapest University of Technology and Economics Institute of Mathematics, \\
Department of Geometry, \\
H-1521 Budapest, Hungary. \\
E-mail:~szirmai@math.bme.hu \\
http://www.math.bme.hu/ $^\sim$szirmai}

\end{document}